\documentclass[preprint, 3p, sort&compress]{elsarticle}

\usepackage{hyperref}
\usepackage{amsmath}
\usepackage{amsfonts}
\usepackage{tikz}

\usetikzlibrary{arrows.meta}
\tikzset{%
	>={Latex[width=2mm,length=2mm]},
	base/.style = {rectangle, rounded corners, draw=black,
		minimum width=2cm, minimum height=.75cm,
		text centered,font=\rmfamily\footnotesize},
	output/.style = {base, fill=blue!15},
	input/.style = {base, fill=green!15},
	process/.style = {base, fill=black!15},
	truth/.style = {base, fill=white}
}

\journal{Reliability Engineering \& System Safety}

\begin{document}

\begin{frontmatter}

\title{Decision-theoretic reliability sensitivity}


\author{Daniel Straub, Max Ehre, Iason Papaioannou}
\address{Engineering Risk Analysis Group, \\ Technische Universit{\"a}t M{\"u}nchen, Germany}

%
%

\begin{abstract}
We propose and discuss sensitivity metrics for reliability analysis, which are based on the value of information. These metrics are easier to interpret than other existing sensitivity metrics in the context of a specific decision and they are applicable to any type of reliability assessment, including those with dependent inputs. We develop computational strategies that enable efficient evaluation of these metrics, in some scenarios without additional runs of the deterministic model. The metrics are investigated by application to numerical examples. 
\end{abstract}

\begin{keyword}
Sensitivity analysis \sep value of information \sep reliability 
\end{keyword}

\end{frontmatter}

\graphicspath{{Figures/}}


\section{Introduction}

The evaluation of sensitivities is an essential part of scientific and engineering analysis. Sensitivities provide information on the relative importance of model input quantities and they support optimization as well as model checking. A large number of sensitivity metrics have been proposed in the literature \cite{saltelli2008global,Iooss2015a,Borgonovo2016}. Metrics exist for deterministic model inputs and for random model inputs. One can distinguish local vs. global metrics. Among the latter, there are variance-based \cite{sobol1993sensitivity,Jansen1999,Saltelli2010} or distribution-based \cite{Chun2000,Borgonovo2007,Borgonovo2016a} metrics. 

The proper choice of a sensitivity metric depends on the application and decision context. In this contribution, we consider sensitivity analysis for random input quantities $\mathbf{X}$ to models that are utilized to evaluate the reliability of a system. Formally, we consider a model $g(\mathbf{X})$, whose output describes the performance of the system in such a way that a value of $g(\mathbf{X})$ below some threshold corresponds to failure $F$. Without loss of generality, this threshold can be set to zero, so that failure is $F=\left\{g(\mathbf{X})\le 0\right\}$ and the probability of failure is $p_F=\Pr\left[g(\mathbf{X})\le 0\right]$. In structural reliability, $g$ is called limit-state function (LSF); in system reliability, $g$ is the structure function.  

Many authors have proposed metrics for such a reliability sensitivity analysis \cite[e.g.,][]{Hohenbichler1986,Madsen1988,Bjerager1989,Ditlevsen1996a,kiureghian2005,Kim2018b,chabridon2018reliability,SARAZIN2021107733}. The sensitivity can be with respect to deterministic input parameters $\theta$, including distribution parameters. In this case, the probability of failure $p_F$ can be interpreted as a function of $\theta$ and a common sensitivity metric is the partial derivative of $p_F$ with respect to $\theta$, $\frac{\partial p_F(\theta)}{\partial \theta_i}$ \cite{Bjerager1989,Ditlevsen1996a,Papaioannou2018}. 

In many instances, the interest is in the sensitivity with respect to the input random variables $\mathbf{X}$. A popular metric in this context are the FORM (First-Order Reliability Method) $\alpha$-factors \cite{Hohenbichler1986,kiureghian2005}, or derivatives thereof \cite{Madsen1988,Ditlevsen1996a,kiureghian2005,Kim2018b}. The FORM $\alpha$-factors  are briefly reviewed in Section \ref{sec:FORM_approx}. 
As an alternative, classical global sensitivity metrics, specifically variance-based and distribution-based metrics, have been adapted for reliability analyses \cite{Cui2010,Li2012,Wei2012,Ehre2020}. In addition, quantile-based sensitivity metrics can be utilized in the context of reliability and risk analysis \cite{Kucherenko2019}.

In many application and decision contexts, the motivation for the sensitivity analysis is to understand on which input random variables one should collect more information to reduce their uncertainty. One refers to this setting as \textit{Factor Prioritization} \cite{saltelli2008global}. To provide a meaningful prioritization, one needs to understand the purpose of the model. In an engineering context, models ultimately serve as decision support. Hence the prioritization should consider how the uncertainty in the input random variables affects the decisions taken based on the model outcome, i.e., by a measure of \textit{decision sensitivity} \cite{Felli1998}. In decision analysis, the \textit{value of information} concept quantifies the effect of reducing uncertainty on the optimality of a decision made based on the model outcome \cite{Raiffa1961,Howard1966}. 

The use of value of information for sensitivity analysis dates back to Felli and Hazen \cite{Felli1998}, who proposed the expected value of partial perfect information\footnote{Felli and Hazen, as well as most other authors, refer to the metric as expected value of perfect information (EVPI). In agreement with the general decision analysis literature we use the term EVPPI for the same metric. The reason for this clarification becomes apparent in Section \ref{sec:VOI}.} (EVPPI) as a sensitivity metric in the context of medical decision making. Around the same time, P\"orn \cite{Porn1997} suggested the use of EVPPI as a component importance measure in Probabilistic Safety Assessments (PSA). More recently, also Borgonovo and Cillo \cite{Borgonovo2017}, Fauriat and Zio \cite{Fauriat2018} and Bj{\o}rnsen et al. \cite{Bjornsen2019a} propose and demonstrate value of information concepts for sensitivity analysis for PSA. In PSA, the system is a deterministic function of binary components with availabilities $p_i$, e.g., through a fault tree model; i.e., the $\mathbf{X}$ are Bernoulli random variables and $g$ is the so-called structure function \cite{rausand2004system}. 


In contrast to these works, we focus on the application of value of information for reliability sensitivity when the reliability is estimated based on a general model with discrete or continuous input random variables. Examples are assessments using physics-based system performance models, such as those found in structural reliability applications \cite{Ditlevsen1996a,Melchers1999}. 
In this context, the value of information concept is applied mainly for optimization of inspection and monitoring schemes \cite[e.g.,][]{Straub2005,Pozzi2011,Thons2013a,Zonta2014,Straub2017,jiang2020optimization,van2020value}. Rather surprisingly, its use as a sensitivity metric has received little attention. 

In this paper, we develop the EVPPI sensitivity metric for two key decision contexts. The first one corresponds to a safety assessment case, in which the decision is between accepting the current system or performing an upgrade that significantly increases reliability. 
The second decision context corresponds to a reliability-based design, in which one or more design parameters can be selected from a typically continuous set of alternatives. The definition of the EVPPI for these two decision contexts, together with algorithms for computing them, are presented in Sections \ref{sec:Safety_assessment} and \ref{sec:Reliability_design}, after a general introduction to the value of information in Section \ref{sec:VOI}. We also demonstrate the relation of the proposed reliability sensitivity metrics to existing sensitivity metrics, in particular the FORM sensitivities and reliability-oriented variance-based sensitivity metrics.

Computation of the EVPPI can be computationally demanding, in particular if the analysis involves advanced computer models. This has been addressed by Oakley and co-authors \cite{Oakley2008,Strong2014,Strong2015}, who present computation strategies for evaluating the EVPPI as well as the expected value of sample information (EVSI) through the use of surrogate modeling strategies. In this paper, we present computational strategies to efficiently evaluate the EVPPI in reliability applications for the two decision contexts. In Section \ref{sec:Safety_assessment}, we show that for the safety assessment case, the EVPPI can be obtained by a mere post-processing of the reliability analysis results without additional runs of the model $g$. 

The proposed sensitivity metrics and their computation are investigated and demonstrated by application to examples in Section \ref{sec:Numerical_examples}.




\section{Value of information}
\label{sec:VOI}

The ultimate goal of an engineering analysis is the recommendation of an optimal action or decision. Example decisions are the selection of a system design, the choice of a management strategy or the decision on whether or not an engineering system is safe to be operated. In line with the literature on decision analysis \cite{Raiffa1961,benjamin1970}, $a$ denotes a decision alternative.

Following expected utility theory, the optimal decision $a_{opt}$ is the one maximizing the expected utility or -- equivalently when the utility is linear with loss -- minimizing the expected loss $L$:
\begin{equation}
a_{opt}=\arg \min_a \text{E}_\mathbf{X}\left[L(\mathbf{X},a)\right].
\end{equation}
$\mathbf{X}$ is the vector of random variables affecting the loss and $\text{E}_\mathbf{X}$ is the mathematical expectation with respect to the probability measure of $\mathbf{X}$. 

The loss function $L$ is problem specific. A frequently used loss function is of the quadratic form,
\begin{equation}
\label{eq:loss_quadratic}
L(\mathbf{X},a)=\left[Y(\mathbf{X})-a\right]^2.
\end{equation}
This loss function occurs when the system design $a$ is optimal if it matches the quantity $Y$ and a deviation of $Y$ from $a$ leads to a loss that increases quadratically with the deviation. $Y$ is a function of the random inputs $\mathbf{X}$. 

In reliability applications, the loss function is typically 
\begin{equation}
\label{eq:loss_reliability}
L(\mathbf{X},a)=c(a)+\mathbb{I}\left[g(\mathbf{X},a)\le 0\right]c_F,
\end{equation}
where $c$ is a function giving the cost of implementing the decision alternative $a$ and $c_F$ is the cost of failure. $\mathbb{I}$ is the indicator function, which results in $1$ if the argument holds and $0$ otherwise. The corresponding expected loss is
\begin{equation}
\label{eq:exp_loss_reliability}
\text{E}_\mathbf{X}\left[L(\mathbf{X},a)\right]=c(a)+\Pr\left[g(\mathbf{X},a)\le 0\right]c_F.
\end{equation}

In some cases it is possible to obtain information on $\mathbf{X}$ prior to making the decision $a$. When data $\mathbf{d}$ is available, it can be utilized to update the probability distribution of $\mathbf{X}$ through Bayesian analysis. An a-posteriori optimal decision can then be found as
\begin{equation}
a_{opt|\mathbf{d}}=\arg \min_a \text{E}_{\mathbf{X}|\mathbf{d}}\left[L(\mathbf{X},a)\right],
\end{equation}
wherein $\text{E}_{\mathbf{X}|\mathbf{d}}$ denotes the expected value with respect to the conditional distribution of $\mathbf{X}$ given $\mathbf{d}$.

By subtracting the corresponding a-posteriori expected loss from the expected loss achieved under the a-priori optimal decision $a_{opt}$, one finds the so-called \textit{conditional value of information}:
\begin{equation}
\label{eq:CVOI_general}
CVOI(\mathbf{d}) = \text{E}_{\mathbf{X}|\mathbf{d}}\left[L(\mathbf{X},a_{opt})\right]-\text{E}_{\mathbf{X}|\mathbf{d}}\left[L(\mathbf{X},a_{opt|\mathbf{d}})\right].
\end{equation}
Since $a_{opt|\mathbf{d}}$ results in the minimum a-posteriori expected loss, it is $\text{E}_{\mathbf{X}|\mathbf{d}}\left[L(\mathbf{X},a_{opt|\mathbf{d}})\right] \le \text{E}_{\mathbf{X}|\mathbf{d}}\left[L(\mathbf{X},a_{opt})\right]  $ and the CVOI cannot be negative. 

Initially, the data $\mathbf{d}$ is not yet available. To understand the potential value of collecting the data, one evaluates the expected value of Eq. \ref{eq:CVOI_general} over the distribution of the data. This is the \textit{(expected) value of information}:
\begin{equation}
\label{eq:VOI_general}
EVOI = \text{E}_{\mathbf{d}} \left[CVOI(\mathbf{d}) \right].
\end{equation}
Since the CVOI cannot be negative, also the EVOI cannot be negative.

Computation of the EVOI is non-trivial in the general case. Algorithms for efficient computation have been proposed in the literature in different application contexts \cite[e.g., ][]{Oakley2008,straub2014,Strong2015,Malings2016}.

In the limit, the decision maker is able to obtain perfect information, such that $\mathbf{X}$ becomes known with certainty. For such a clairvoyant decision maker, finding the optimal $a$ reduces to a deterministic decision problem:
\begin{equation}
a_{opt|\mathbf{x}}=\arg \min_a L(\mathbf{x},a).
\end{equation}
By analogy with Eq. \ref{eq:VOI_general}, the \textit{expected value of perfect information} is
\begin{equation}
\label{eq:EVPI_general}
EVPI = \text{E}_{\mathbf{X}} \left[L(\mathbf{X},a_{opt})-L(\mathbf{X},a_{opt|\mathbf{X}}) \right].
\end{equation}
The EVPI is the upper bound on the EVOI, i.e., the best data collection cannot provide a higher expected value than EVPI in the considered decision context. 

A special case of value of information arises when obtaining perfect information on a subset of input parameters. Consider that input random variable $X_i$ is observed to take value $x_i$. One can evaluate the \textit{conditional value of partial perfect information} on $X_i$ as
\begin{equation} \label{eq:CVPPI}
CVPPI_{X_i}(x_i)=\text{E}_{\mathbf{X}_{-i}}\left[L(\mathbf{X},a_{opt})|X_i=x_i\right]-\text{E}_{\mathbf{X}_{-i}}\left[L(\mathbf{X},a_{opt|x_i})|X_i=x_i\right],
\end{equation}
wherein ${\mathbf{X}_{-i}}$ indicates the vector of all random variables $\mathbf{X}$ except $X_i$.

By evaluating the expected value of the CVPPI over the prior distribution of $X_i$, one obtains the expected value of partial perfect information:
\begin{equation}
\label{eq:VPPI}
\begin{split}
EVPPI_{X_i}&= \text{E}_{X_i}\left[CVPPI_{X_i}(X_i)\right] \\ 
&= \text{E}_{\mathbf{X}}\left[L(\mathbf{X},a_{opt})\right] -  \text{E}_\mathbf{X}\left[L(\mathbf{X},a_{opt|X_i})\right].
\end{split}
\end{equation}

The EVPPI is a decision-theoretic metric for factor prioritization. It describes the expected gains due to improved decision making when learning a specific input random variable. 

\textit{Remark 2.1}: When the loss function is of the quadratic form of Eq. \ref{eq:loss_quadratic}, the optimal decision a-priori is
$a_{opt}= \text{E}_{\mathbf{X}}\left[Y(\mathbf{X})\right]$. The associated expected loss is
$ \text{E}_{\mathbf{X}}\left[L(\mathbf{X},a_{opt})\right]=\text{Var}_{\mathbf{X}}\left[Y(\mathbf{X})\right]$. Correspondingly, the expected loss associated with the optimal posterior decision is
$\text{E}_{\mathbf{X}_{-i}}\left[L(\mathbf{X},a_{opt|x_i})|X_i=x_i\right]=\text{Var}_{\mathbf{X}_{-i}}\left[Y(\mathbf{X})|X_i=x_i\right]$. 
Inserting these expressions into Eq. \ref{eq:VPPI} gives
\begin{equation}
\begin{split}
EVPPI_{X_i}&= \text{Var}_{\mathbf{X}}\left[Y(\mathbf{X})\right] -  \text{E}_{X_i} \left\{\text{Var}_{\mathbf{X}_{-i}}\left[Y(\mathbf{X})|X_i\right]\right\} \\
&= \text{Var}_{X_i} \left\{\text{E}_{\mathbf{X}_{-i}}\left[Y(\mathbf{X})|X_i\right]\right\}.
\end{split}
\end{equation}
Dividing this result with $\text{Var}_{\mathbf{X}}\left[Y(\mathbf{X})\right]$ gives the first-order variance-based sensitivity index, the Sobol' index\footnote{Under a quadratic loss function, $\text{Var}_{\mathbf{X}}\left[Y(\mathbf{X})\right]$ is the expected value of perfect information EVPI according to Eq. \ref{eq:EVPI_general}, since the expected loss under perfect information is zero. Therefore, the first-order Sobol' index corresponds to the EVPPI normalized with the EVPI.}. Hence, the EVPPI associated with a quadratic loss function is equivalent to the first-order Sobol' index \cite{Oakley2008}.
The EVPPI is also related to other sensitivity metrics, as discussed in \cite{Borgonovo2016,borgonovo2021probabilistic}.

\textit{Remark 2.2}: The definition of the EVPPI can be extended to groups of input random variables $\mathbf{X}_{\boldsymbol{v}}=\left\{X_i,\,i\in \boldsymbol{v} \right\}$, with $\boldsymbol{v}$ being a subset of $\left\{1,\,\dots,\,n\right\}$. In the case of the quadratic loss function, the resulting $EVPPI_{\mathbf{X}_{\boldsymbol{v}}}$ normalized with $\text{Var}_{\mathbf{X}}\left[Y(\mathbf{X})\right]$ is equal to the closed Sobol' index of $\mathbf{X}_{\boldsymbol{v}}$.

The formulation of the decision alternatives and the loss function depends on the decision context. In Section \ref{sec:Safety_assessment}, we derive the EVPPI for cases when the analysis is performed to assess the safety of a given system or design. In Section \ref{sec:Reliability_design} we consider the case of reliability-based optimization, in which $a$ includes a set of optimization parameters. For both cases, we discuss efficient estimators for the EVPPI. In the context of reliability sensitivity, efficiency is measured in terms of the number of additional evaluations of $g(\mathbf{x})$ after the reliability is computed. 

\section{Safety assessment}
\label{sec:Safety_assessment}
\subsection{Decision analysis}
\begin{figure}
	\centering	\includegraphics[width=80mm]{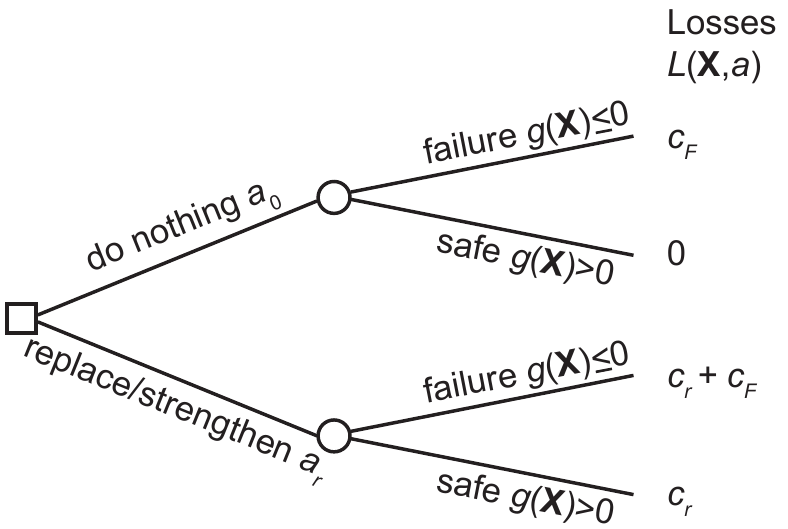}
	\caption{Decision tree representing the safety assessment.}
	\label{fig:decision_tree_safety}
\end{figure}

A common decision situation in engineering is associated to the safety assessment of an existing system. In its simplest form, the problem can be represented by the decision tree of Figure \ref{fig:decision_tree_safety}. It has two decision alternatives $a=\{a_0,a_r\}$, corresponding to doing nothing and replacing/strengthening the system. The state of the system is binary, it either fails $F=\left\{g(\mathbf{X})\le 0 \right\}$ or survives $\bar{F}=\left\{g(\mathbf{X})> 0 \right\}$.
The loss function is of the form given by Eq. \ref{eq:loss_reliability}.

We make the assumption that following a replacement or strengthening action $a_r$, the system has a negligible probability of failure\footnote{While this is not the case in reality, it can be expected that the probability of failure of the updated system $\Pr(F|a_r)$ is indeed much lower than the original $\Pr(F)$.}, hence $\Pr(F|a_r)=0$. This model is equivalent to the one of \cite{Porn1997}, in which $a_r$ is the decision to reject a component in a system design. The resulting expected loss is
\begin{equation} \label{eq:safety_exp_loss}
\text{E}_{\mathbf{X}}\left[L(\mathbf{X},a)\right]= 
\begin{cases}
p_F c_F, & a=a_0 \\
c_r & a=a_r,
\end{cases}       
\end{equation}
where $p_F=\Pr[g(\mathbf{X})\le 0]$, $c_F$ is the cost of failure and $c_r$ is the cost of replacement.

In this case, the optimal decision a-priori is
\begin{equation}
a_{opt}= 
\begin{cases}
a_0, & p_F\leq \frac{c_r}{c_F} \\
a_r, & \text{else},
\end{cases}       
\end{equation}
wherein $\frac{c_r}{c_F}$ serves as a decision threshold.

When $X_i$ is known, the analysis is performed with $p_{F}(x_i)=\Pr\left[g(\mathbf{X})\le 0|X_i=x_i\right]$, the conditional probability of failure given $X_i=x_i$. The conditionally optimal decision is
\begin{equation}
a_{opt|X_i}(x_i)= 
\begin{cases}
a_0, & p_F(x_i)\leq \frac{c_r}{c_F} \\
a_r, & \text{else}.
\end{cases}       
\end{equation}

If the conditional $p_{F}(x_i)$ is on the same side of the threshold $\frac{c_r}{c_F}$ as the unconditional $p_{F}$, the optimal decision is not altered by the knowledge $X_i=x_i$. In this case, the CVPPI is zero. The CVPPI is only positive, if the decision is changed. Hence it is
\begin{equation} \label{eq:safety_CVPPI}
CVPPI_{X_i}(x_i)= 
\begin{cases} \left|c_F p_{F}(x_i)-c_r\right|, & a_{opt|X_i}(x_i)\ne a_{opt} \\
0, &a_{opt|X_i}(x_i) = a_{opt}.
\end{cases}       
\end{equation}

\begin{figure}
	\centering	\includegraphics[width=80mm]{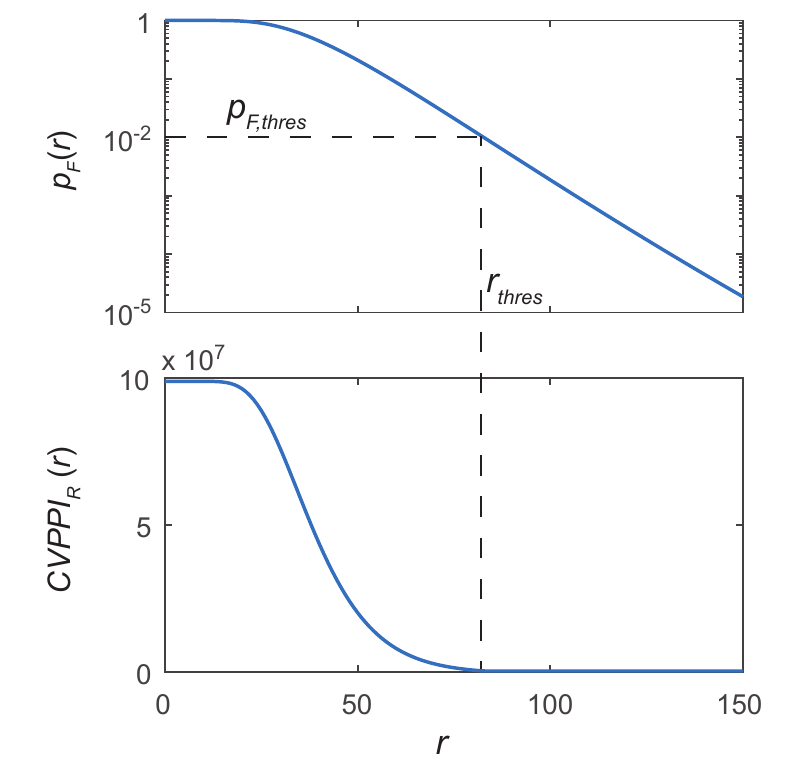}
	\caption{Conditional probability of failure for input random variable $R$ of the example in Section \ref{sec:Example_component_problem}, together with the conditional value of partial perfect information CVPPI. The cost of failure is $c_F=10^8$ and the cost of replacement $c_r=10^6$, hence the threshold is $\frac{c_r}{c_F}=10^{-2}$. The prior probability of failure is $\Pr(F)=7.4\cdot10^{-3}$, hence $a_{opt}=a_0$ and the CVPPI is non-zero only for values of  $p_F(r)>10^{-2}$.} 
	\label{fig:CVPPI_linear_Example}
\end{figure}

The CVPPI is illustrated in Figure \ref{fig:CVPPI_linear_Example}.

The resulting EVPPI is
\begin{equation} \label{eq:VPPI_safety_1}
EVPPI_{X_i}= 
\int_{\Omega_{X_i}}
\left|c_F p_{F}(x_i)-c_r\right|f_{X_i}(x_i) \, dx_i,
\end{equation}
wherein the integration domain $\Omega_{X_i}$ is the set of values $x_i$ for which the conditionally optimal decision $a_{opt|X_i}(x_i)$ differs from the unconditional optimum $a_{opt}$.

For example, if the optimal prior decision is $a_{opt}=a_0$, then the domain $\Omega_{X_i}$ consists of all $x_i$ for which the following inequality holds:
\begin{equation}
 p_F(x_{i}) > \frac{c_r}{c_F}.
\end{equation}

In most cases, $p_F(x_i)$ is a monotonic function. In these cases, a threshold value $x_{i,thres}$ can be found at which the optimal decision changes from $a_0$ to $a_r$, or reversely. It can be determined from the condition
\begin{equation}
 p_F(x_{i,thres}) = \frac{c_r}{c_F}.
\end{equation}

The evaluation of the EVPPI following Eq. \ref{eq:VPPI_safety_1} requires the function $p_F(x_i)$. Two strategies for estimating this function are presented in the following subsections.

\textit{Remark 3.1}: When evaluating the EVPPI for a group of input random variables $\mathbf{X}_{\boldsymbol{v}}=\left\{X_i,\,i\in \boldsymbol{v} \right\}$,
Eqs. \ref{eq:safety_exp_loss} to \ref{eq:safety_CVPPI} hold with $x_i$ replaced by $\mathbf{x}_{\boldsymbol{v}}$. The EVPPI is 
\begin{equation} \label{eq:VPPI_safety_group}
EVPPI_{\mathbf{X}_{\boldsymbol{v}}}= 
\int_{\Omega_{\mathbf{X}_{\boldsymbol{v}}}}
\left|c_F p_{F}(\mathbf{x}_{\boldsymbol{v}})-c_r\right|f_{\mathbf{X}_{\boldsymbol{v}}}(\mathbf{x}_{\boldsymbol{v}}) \, d\mathbf{x}_{\boldsymbol{v}},
\end{equation}
wherein the integration domain $\Omega_{\mathbf{X}_{\boldsymbol{v}}}$ is the set of values for which $a_{opt|\mathbf{X}_{\boldsymbol{v}}}(\mathbf{x}_{\boldsymbol{v}})\ne a_{opt}$.

The computational strategies for evaluating $EVPPI_{X_i}$ that we present in the following can be extended to the computation of $EVPPI_{\mathbf{X}_{\boldsymbol{v}}}$, but we do not discuss this further.

\subsection{Estimation based on failure samples} \label{sec:failure_samples}
When estimating $p_F$ with a sampling-based method, the function $p_F(x_i)$ can be estimated from samples in the failure domain \cite{Li2019a}. To this end, we note that the conditional $p_{F}(x_i)=\Pr\left[F|X_i=x_i\right]$ can be found through Bayes' rule as
\begin{equation} 
p_F(x_i)=\frac{f_{X_i|F}(x_i) }{f_{X_i}(x_i)} p_F,
\label{eq:pF_cond_Bayes}
\end{equation}
wherein $f_{X_i|F}$ is the PDF of $X_i$ conditional on a failure event. 

\begin{figure}
	\centering	\includegraphics[width=80mm]{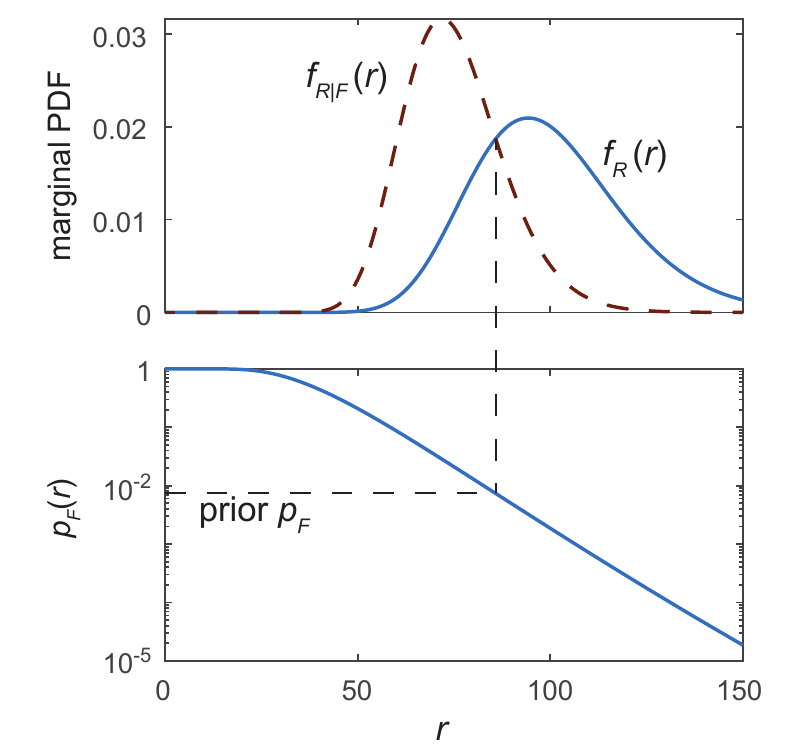}
	\caption{Conditional and unconditional PDF of the input random variable $R$ of the example in Section \ref{sec:Example_component_problem}, together with the conditional failure probability $p_F(r)=\Pr(F|R=r)$.}
	\label{fig:ConditionalPDF_Example_1}
\end{figure}

The relationship of Eq. \ref{eq:pF_cond_Bayes} is illustrated in Figure \ref{fig:ConditionalPDF_Example_1}. The fact that the EVPPI is a function of the relation between $f_{X_i|F}$ and $f_{X_i}$ points to the link to distribution-based sensitivity metrics, which is discussed in \cite{Borgonovo2016c}.

\begin{figure}
	\centering	\includegraphics[width=80mm]{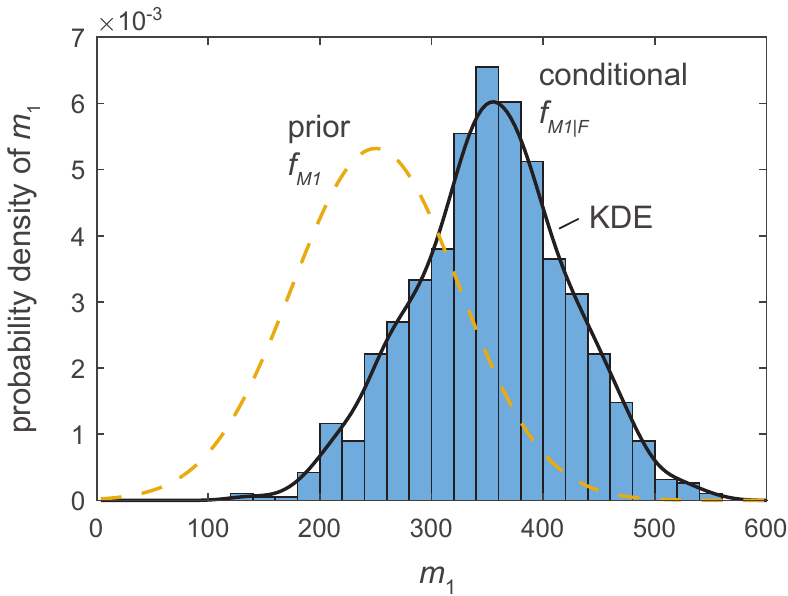}
	\caption{Histogram of failure samples of the input random variable $M_1$ obtained as a side product from a Monte-Carlo-based reliability analysis, together with the corresponding kernel density estimate KDE. Taken from the example in Section \ref{sec:nonlinear_example} The number of failure samples is $n_F=947$. Additionally, the unconditional (prior) PDF of $M_1$ is shown.}
	\label{fig:KDE_fit_illustration}
\end{figure}

If crude Monte Carlo simulation is utilized to estimate $p_F$, then the set of samples that fall into the failure domain are independent and identically distributed (iid) samples from $f_{X_i|F}$. Hence they can be utilized to obtain an estimate of $f_{X_i|F}$, e.g., by means of a kernel density estimator (KDE) as illustrated in Figure \ref{fig:KDE_fit_illustration}.

If more advanced sampling techniques are employed to estimate $p_F$, an additional step might be necessary before estimating $f_{X_i|F}$. In the case of importance sampling techniques, the resulting samples in the failure domain are weighted samples. Approximate iid samples from $f_{X_i|F}$ can then be obtained by an additional resampling step. When employing subset simulation, the resulting samples in the failure domain are correlated. In this case, the same procedure as with crude Monte Carlo is applicable, but the quality of the estimates will be lower with the same number of samples. The estimation of $f_{X_i|F}$ from failure samples obtained with different sampling techniques is presented in \cite{Li2019a}. 


\subsection{FORM-based approximation}
\label{sec:FORM_approx}
A powerful approximation method to compute $\Pr[F]$ is the First-order Reliability Method (FORM) \cite{Rackwitz2001b}. 
FORM is based on a transformation of the limit-state function  $g(\mathbf{X})$ to an equivalent $G(\mathbf{U})$ defined in terms of independent standard normal random variables $\mathbf{U}$. This is achieved by an iso-probabilistic transform, typically the Rosenblatt \cite{Hohenbichler1981a} or the Nataf \cite{Liu1986} transform. The probability of failure is then $p_F=\Pr\left[G\left(\mathbf{U}\right)\le 0\right]$ and can be approximated by a linearization of  $G\left(\mathbf{U}\right)$ at the so-called most likely failure point $\mathbf{u}^*$. This is the point in the failure domain $\left\{G\left(\mathbf{U}\right)\le 0\right\}$ with the highest probability density. Because of the rotational symmetry of the standard normal PDF, this is also the point in $\left\{G\left(\mathbf{U}\right)\le 0\right\}$ closest to the origin:
\begin{equation}
\begin{split}
	\mathbf{u}^*=& \arg\min_\mathbf{u} \lVert \mathbf{u} \rVert \\
    & \text{s.t. } G\left(\mathbf{u}\right)\le 0.
\end{split}
\end{equation}
$ \lVert \cdot \rVert$ is the Euclidean norm.

Let $G_1$ denote the linearized limit-state function. Provided that $\mathbf{u}^*$ is not the origin, it can be shown that 
$\Pr\left[G_1\left(\mathbf{U}\right)\le 0\right]=\Phi\left[-\beta_0\right]$, where $\beta_0 = \lVert \mathbf{u}^* \rVert$ is the distance of the most likely failure point from the origin and $\Phi$ is the standard normal CDF \cite{Ditlevsen1996a,kiureghian2005}. Hence the FORM approximation of the probability of failure is
\begin{equation}
	p_F \approx \Phi\left[-\beta_0\right].
\end{equation}

The FORM $\alpha$-factors $\boldsymbol{\alpha}=[\alpha_1;\dots;\alpha_n]$ are the directional cosines of $\mathbf{u}^*$ and are commonly used as a sensitivity metric \cite{Hohenbichler1986,kiureghian2005,PapaioannouFORMsensitivity}. 

To evaluate the EVPPI, we note that the FORM estimate of the probability of failure conditional on $U_i=u_i$ is
\begin{equation}\label{eq:pF_ui}
\begin{split}
p_F(u_i)
&=\Phi\left(-\frac{\beta_0-\alpha_i u_i}{\sqrt{\alpha_1^2+\alpha_2^2+\dots+\alpha_{i-1}^2+\alpha_{i+1}^2+\dots+\alpha_n^2}}\right) \\ 
&= \Phi\left(\frac{\alpha_i u_i-\beta_0}{\sqrt{1-\alpha_i^2}}\right). 
\end{split} 
\end{equation}

In the case of independent random variables, $X_i$ is related to $U_i$ by the marginal transformation 
\begin{equation}
u_i = \Phi^{-1}[F_{X_i}(x_i)],
\end{equation}
hence in this case it is
\begin{equation}
\label{eq:FORM_p_F_x_i}
p_F(x_i)=\Phi\left\{\frac{\alpha_i \Phi^{-1}[F_{X_i}(x_i)]-\beta_0}{\sqrt{1-\alpha_i^2}}\right\}.
\end{equation}
Setting $p_F(x_{i,thres})=\frac{c_r}{c_F}$ results in
\begin{equation} \label{eq:FORM_x_i_thresh}
x_{i,thres}=F_{X_i}^{-1}\left(\Phi\left\{\frac{1}{\alpha_i}\left[\sqrt{1-\alpha_i^2} \Phi^{-1}\left(\frac{c_r}{c_F}\right)+\beta_0\right]\right\}\right).
\end{equation}

The threshold of Eq. \ref{eq:FORM_x_i_thresh} is the boundary of the integration domain in the EVPPI calculation,  Eq. \ref{eq:VPPI_safety_1}. Inserting Eq. \ref{eq:FORM_p_F_x_i} into Eq. \ref{eq:VPPI_safety_1} results in a FORM estimate of the EVPPI. Integration by substitution with $\Phi^{-1}\left[F_{X_i}(x_i)\right]=u_i$ gives
\begin{equation} \label{eq:EVPPI_safety_FORM}
EVPPI_{X_i}= 
\int_{\Omega_{u_i}}
\left|c_F p_{F}(u_i)-c_r\right|\varphi\left(u_i\right) \, du_i,
\end{equation}
wherein $p_{F}(u_i)$ is given by Eq. \ref{eq:pF_ui}, $\varphi$ is the standard normal PDF and the integration domain is
\begin{equation}
\Omega_{u_i}=
\begin{cases}
\left[u_{i,thres},\infty\right], & p_F\leq \frac{c_r}{c_F}, \, \alpha_i>0 
\\
\left[-\infty,u_{i,thres}\right], & p_F\leq \frac{c_r}{c_F}, \, \alpha_i<0 
\\
\left[-\infty,u_{i,thres}\right], & p_F> \frac{c_r}{c_F} , \, \alpha_i>0 \\
\left[u_{i,thres},\infty\right], & p_F> \frac{c_r}{c_F} , \, \alpha_i<0.
\end{cases}
\end{equation}
$u_{i,thres}$ is the transformation of the threshold $x_{i,thresh}$ to standard normal space: 
\begin{equation} \label{eq:v_thresh}
	u_{i,thres}=\frac{1}{\alpha_i}\left[\sqrt{1-\alpha_i^2} \Phi^{-1}\left(\frac{c_r}{c_F}\right)+\beta_0\right].
\end{equation}

In \cite{Papaioannou2022}, we show that the integral of Eq. \ref{eq:EVPPI_safety_FORM} can alternatively be evaluated as 
\begin{equation}
EVPPI_{X_i}=
\lvert c_F \Phi_2\left(-\beta_0,s_i u_{i,thres},-s_i\alpha_i\right)-c_r\Phi\left(s_i u_{i,thres}\right)\rvert.
\end{equation}
$\Phi_2(x_1,x_2,r)$ is the bivariate standard normal CDF with correlation coefficient $r$ evaluated at $x_1$ and $x_2$ and $s_i$ is either $=+1$ or $-1$, depending on the integration domain: $s_i=sgn\left[\left(p_F-\frac{c_r}{c_F}\right)\alpha_i\right]$ if $\left(p_F-\frac{c_r}{c_F}\right)\alpha_i\ne 0$, otherwise $s_i=-1$.

\textit{Remark 3.2}: The FORM approximation of the EVPPI is identical for $+\alpha_i$ and $-\alpha_i$, hence it is also possible to utilize $|\alpha_i|$ in the above expressions.


\subsection{Relation of EVPPI to the FORM sensitivity index}

Equations \ref{eq:EVPPI_safety_FORM}--\ref{eq:v_thresh} show that the FORM approximation of the EVPPI does not depend on the marginal distribution of $X_i$. It is a function only of $\beta_0$, $\alpha_i$ and $\frac{c_r}{c_F}$. Hence, one can directly obtain the EVPPI corresponding to a FORM sensitivity $\alpha_i$ without any calls of the $g$-function.

\begin{figure}
	\centering	\makebox[0pt]{\includegraphics[width=80mm]{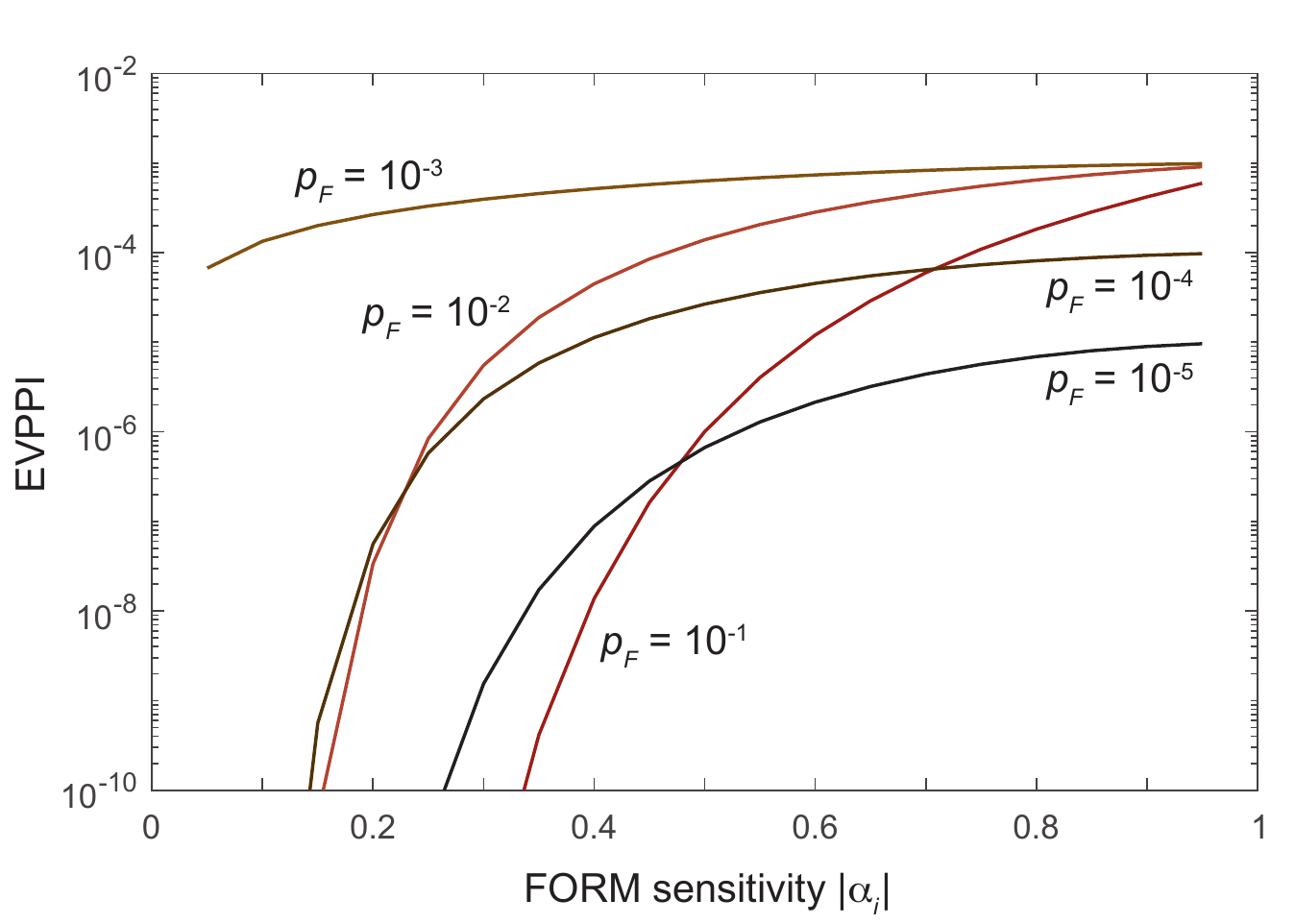}}
	\caption{EVPPI in function of the FORM sensitivity index $\alpha_i$, for varying values of $p_F$, $c_F=1$ and $\frac{c_r}{c_F}=10^{-3}$.}
	\label{fig:FORM_VPPI_cr10e3}
\end{figure}

\begin{figure}
	\centering	\makebox[0pt]{\includegraphics[width=80mm]{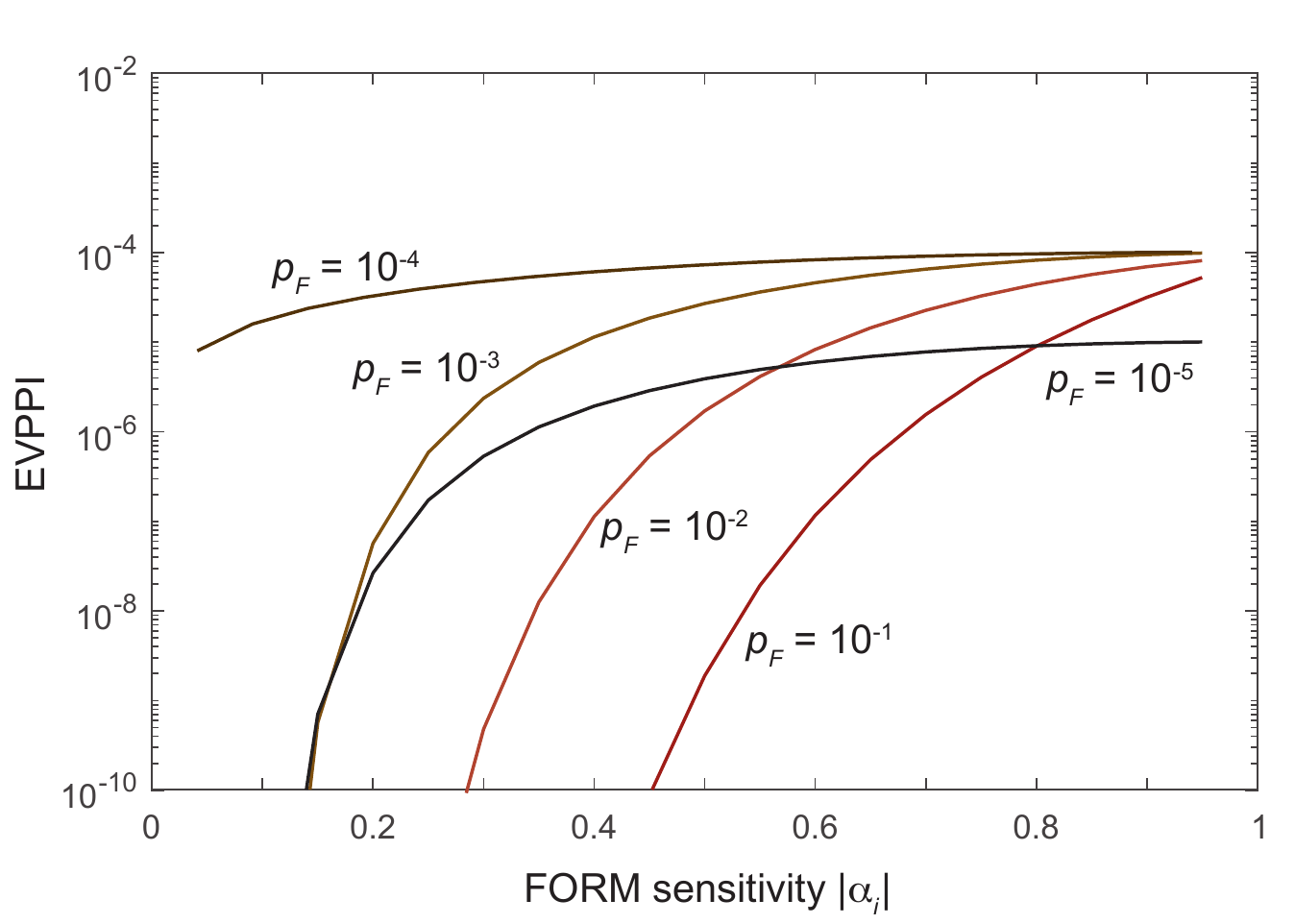}}
	\caption{EVPPI in function of the FORM sensitivity index $\alpha_i$, for varying values of $p_F$, $c_F=1$ and $\frac{c_r}{c_F}=10^{-4}$.}
	\label{fig:FORM_VPPI_cr10e4}
\end{figure}

Figure \ref{fig:FORM_VPPI_cr10e3} shows the EVPPI in function of $|\alpha_i|$ for different values of $p_F$, with cost ratio $\frac{c_r}{c_F}=10^{-3}$; Figure \ref{fig:FORM_VPPI_cr10e4} shows the same for cost ratio $\frac{c_r}{c_F}=10^{-4}$. The EVPPI is highest when $p_F$ is close to $\frac{c_r}{c_F}$. This is to be expected, because in these cases the observation $X_i=x_i$ is most likely to lead to a change in the optimal decision.

The EVPPI is increasing with increasing value of $|\alpha|$. 
Since $\beta_0$ and $\frac{c_r}{c_F}$ are fixed for a given reliability problem, the ordering according to FORM approximation of the EVPPI is the same as the ordering according to the absolute values of the FORM sensitivities $|\alpha|$. This holds for reliability problems with independent input random variables when the decision problem is as in Figure \ref{fig:decision_tree_safety}. 

The key benefit of using the EVPPI metric rather than the FORM $\alpha$-factors is that it accounts for the effect of $\beta_0$ and $\frac{c_r}{c_F}$ on the decision sensitivity. For example, consider the case $\frac{c_r}{c_F}=10^{-3}$  shown in Figure \ref{fig:FORM_VPPI_cr10e3}. One can observe that for $p_F=10^{-3}$ an input random variable with $|\alpha|=0.8$ has an EVPPI that is double the one associated with $|\alpha|=0.35$. In contrast, for $p_F=10^{-2}$ the EVPPI corresponding to $|\alpha|=0.8$ is $34$ times larger than the EVPPI corresponding to $|\alpha|=0.35$. Hence the ratio between $|\alpha|$ (or $\alpha^2$) values is not a good indicator for how much more important one input random variable is relative to another in the specific decision context. But the ratio between EVPPI values is; it is admissible to make a statement like ``learning $X_i$ is $\frac{EVPPI_{X_i}}{EVPPI_{X_j}}$ more valuable than learning $X_j$''.

\subsection{Normalization}\label{sec:normalization_safety}
It can be preferable to express the sensitivity by an index that is normalized to lie between $0$ and $1$. A trivial normalization can be achieved by normalizing the EVPPI by the sum of the EVPPI of all input random variables:
\begin{equation}\label{eq:normalized_EVPPI}
\text{Normalized }EVPPI_{X_i}= \frac{EVPPI_{X_i}}{\sum_{j=1}^n EVPPI_{X_j}}.
\end{equation}
The advantage of this normalization is that the indices sum up to one. Therefore, this normalization is suitable for comparison with other sensitivity indices that have this property or are normalized accordingly. 

A natural normalization is obtained by evaluating the EVPPI relative to the expected value of perfect information (EVPI) of Eq. \ref{eq:EVPI_general}. The EVPI is the value associated with having full information on all inputs. 
We refer to this normalized index as the \textit{relative} EVPPI:
\begin{equation}\label{eq:relative_EVPPI}
\text{Relative }EVPPI_{X_i}= \frac{EVPPI_{X_i}}{EVPI}.
\end{equation}

For the safety assessment decision, the EVPI is
\begin{equation} \label{eq:safety_EVPI}
EVPI= 
\begin{cases}
p_F \left(c_F-c_r \right), & p_F\leq \frac{c_r}{c_F} \\
c_r \left(1-p_F \right), & \text{else},
\end{cases}           
\end{equation}
and is readily available.

Because the EVPI is an upper bound on the EVPPI, the relative EVPPI is bounded by zero and one. However, the relative EVPPIs do not sum up to one. Furthermore, their sum can be larger than one, because the EVPPI for groups of input random variables can be lower than the sum of the EVPPIs of individual input random variables \cite{brennan2007calculating}. 

\textit{Remark 3.3}: As discussed in Section \ref{sec:VOI}, the Sobol' index is a relative EVPPI for the case of the quadratic loss function.  

\section{Reliability-based design}
\label{sec:Reliability_design}
\subsection{Decision analysis}
In reliability-based design, the decision is typically not binary as in the case of the safety assessment. Instead, one or more design parameters $\mathbf{a}$ can be selected from an (often continuous or continuous-discrete) domain. 
The function $c_d(\mathbf{a})$ describes the design cost, which is here assumed to be deterministic, i.e., it does not depend on $\mathbf{X}$. As in Section \ref{sec:Safety_assessment}, the cost of a failure is $c_F$. Hence the loss function is
\begin{equation}
L(\mathbf{X},\mathbf{a})=c_d(\mathbf{a})
+\mathbb{I}\left[g(\mathbf{X},\mathbf{a})\le 0\right]c_F.
\end{equation}

The optimal decision a-priori is
\begin{equation} \label{eq:prior_optimization_design}
\mathbf{a}_{opt}
= \arg \min_\mathbf{a} \, c_d(\mathbf{a}) + p_F(\mathbf{a}) c_F .
\end{equation}

The a-posteriori optimal decision given $X_i=x_i$ is
\begin{equation} \label{eq:posterior_optimization}
\mathbf{a}_{opt|X_i}(x_i)= \arg \min_\mathbf{a} \, c_d(\mathbf{a}) + p_F(x_i,\mathbf{a}) c_F .
\end{equation}
$p_F(x_i,\mathbf{a})=\Pr\left[g(\mathbf{X},\mathbf{a})\le 0|X_i=x_i\right]$ is the conditional probability of failure given $X_i=x_i$ for a design $\mathbf{a}$. 

Following Eq. \ref{eq:VPPI}, the EVPPI is
\begin{equation}
\label{eq:VPPI_design}
EVPPI_{X_i}= c_d(\mathbf{a}_{opt}) + p_F(\mathbf{a}_{opt}) c_F   -  \text{E}_{X_i}\left\{c_d[\mathbf{a}_{opt|X_i}(X_i)] + p_F[X_i,\mathbf{a}_{opt|X_i}(X_i)] c_F\right\}.
\end{equation}

This computation requires that $p_F(x_i,\mathbf{a})$ is available in a manner that facilitates the solution of Eq. \ref{eq:posterior_optimization}. This optimization problem has to be solved many times for varying values of $X_i$, hence the evaluation of $p_F(x_i,\mathbf{a})$ must be efficient and without noise. In general, $p_F(x_i,\mathbf{a})$ will not be available in analytical form. We present a computational strategy to address this challenge in the following. In Section \ref{sec:FORM_design_case}. we additionally discuss a FORM approximation.

\subsection{Discrete or discretized design choices}\label{sec:discretized_design}

In the case where the domain of $\mathbf{a}$ is discrete, or can be discretized, one can solve $m$ reliability problems, one for each design choice $\{a_1,\,\dots \,a_m\}$. For a given $a_j$, the conditional $p_F(x_i,a_j)$ can be determined with the strategies presented in Sections \ref{sec:failure_samples} and \ref{sec:FORM_approx}. 
The solution of Eq. \ref{eq:posterior_optimization} then reduces to selecting the optimal $a_j$ from the discrete set $\{a_1,\,\dots \,,a_m\}$.

\begin{figure}
	\centering	\makebox[0pt]{\includegraphics[width=80mm]{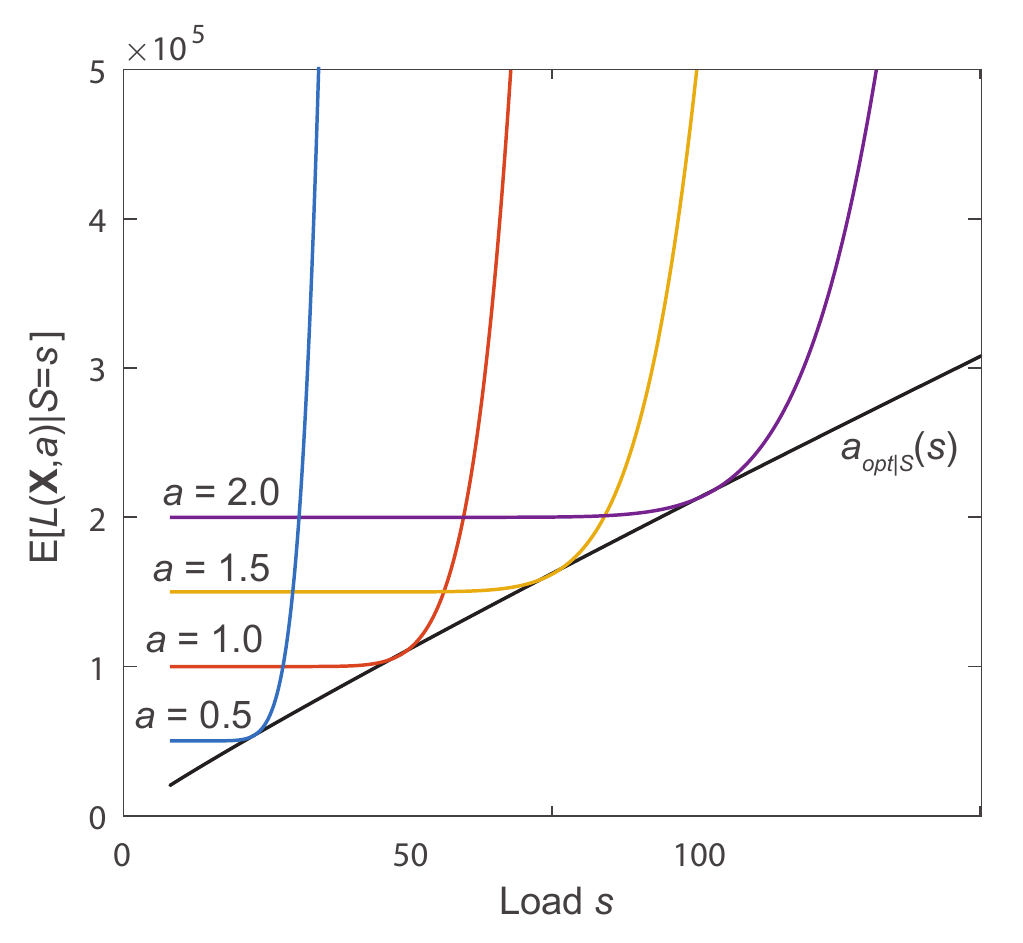}}
	\caption{Expected loss $\text{E}\left[L(\mathbf{X},a)|S=s\right]$ conditional on input random variable $S$, from the example in Section \ref{sec:Example_component_problem}. The black line is the posterior expected loss under the a-posteriori optimal choice, $\text{E}\left[L(\mathbf{X},a_{opt|s}(s))|S=s\right]$. In contrast, the colored lines represent the posterior expected loss for fixed design choices, $\text{E}\left[L(\mathbf{X},a_j)|S=s\right]$. The minimum of these would apply for a specific value of $s$. The black and the colored lines meet at the value of $s$ at which the fixed choice $a_i$ is identical to the a-posteriori optimal choice $a_{opt|s}$.}
	\label{fig:Posterior_loss_discrete_Ex1}
\end{figure}

When the domain of $\mathbf{a}$ is continuous, the consideration of only a discrete set of design choices $\{a_1,\,\dots \,a_m\}$ leads to an approximation error. The difference between the exact and approximated expected a-posteriori loss $\text{E}_{\mathbf{X}_{-i}}[L(\mathbf{X},a_{opt|X_i}|{X_i=x_i})]=c_d[\mathbf{a}_{opt|X_i}(x_i)] + p_F[x_i,\mathbf{a}_{opt|X_i}(x_i)] c_F$, is illustrated in Figure \ref{fig:Posterior_loss_discrete_Ex1}. The expected loss under the a-posteriori optimal choice $a_{opt|X_i}$ is overestimated for all values of $X_i$, except for those where the discrete choices $a_j$ coincide with the optimal choice in the continuous domain of $\mathbf{a}$.

To limit the error from this approximation, the discrete design choices should cover all values of $\mathbf{a}$ that might potentially be optimal. The spacing between $a_i$'s should represent a good trade-off between accuracy and computational cost. These choices are investigated numerically in Section \ref{sec:example_discrete_approx}.  

\subsection{A remark on the FORM approximation}
\label{sec:FORM_design_case}
It is possible to derive a FORM approximation of the EVPPI also in the design decision case. Such an approximation is described in the Annex of the paper, assuming a linear design LSF.

\begin{figure}
	\centering	\makebox[0pt]{\includegraphics[width=80mm]{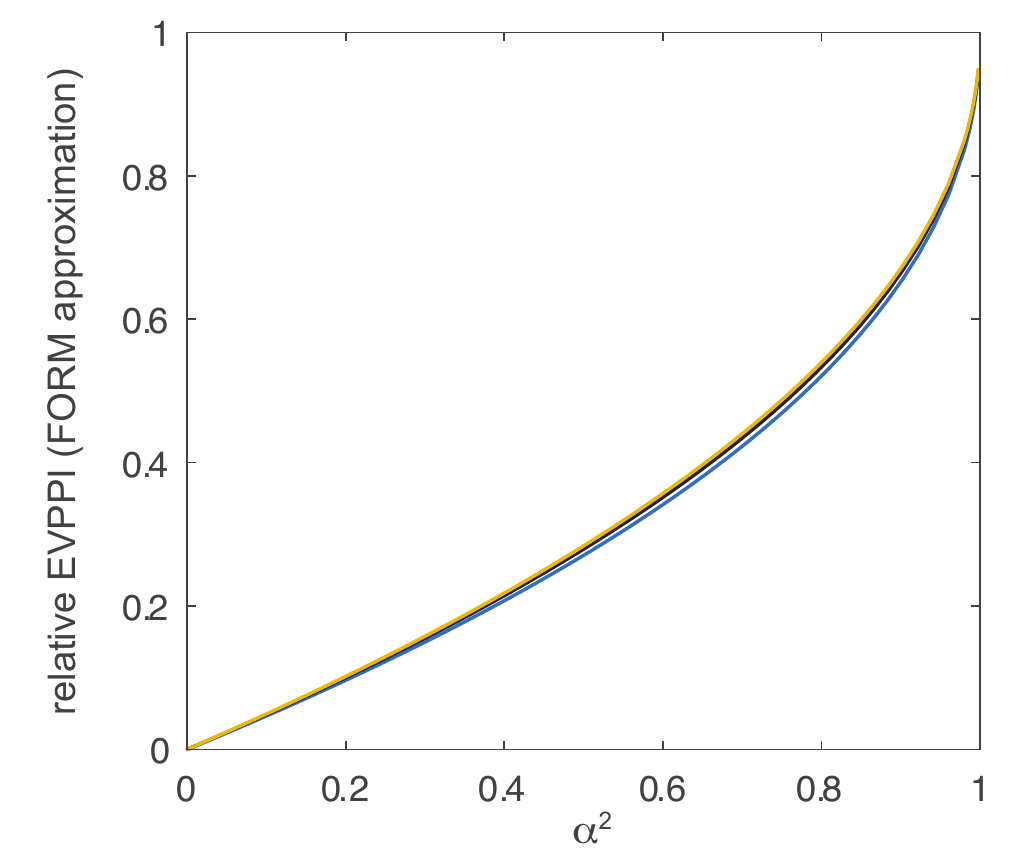}}
	\caption{EVPPI for the reliability-based design case evaluated from the FORM approximation, normalized with the EVPPI at $\alpha^2=1$. The three curves correspond to a reliability index $\beta_0=3$ (blue), $4$ (black) and $5$ (yellow). Details are provided in the Annex.}
	\label{fig:FORM_design_norm}
\end{figure}

The EVPPI for different values of the squared FORM $\alpha$-values relative to the EVPPI of $\alpha^2=1$ are shown in Figure \ref{fig:FORM_design_norm}. These results are independent of the cost parameters $c_\delta$ and $c_F$ introduced in the Annex.

Based on the results of the Annex and Figure \ref{fig:FORM_design_norm} one can conclude that for any design LSF with independent inputs, whose reliability is well approximated by FORM, one should expect the relative EVPPI of the design decision case to give similar importance measures than the FORM $\alpha$-factors. 

\subsection{Normalization}
By analogy with Section \ref{sec:normalization_safety}, the normalized EVPPI and the relative EVPPI can also be defined for the reliability-based design. While computation of the former is trivial, the computation of the EVPI necessary for determining the relative EVPPI generally requires additional model runs. In most cases, this additional effort is not justified.

\section{Numerical investigations}
\label{sec:Numerical_examples}

\subsection{Example 1: Component reliability}
\label{sec:Example_component_problem}

\subsubsection{VOI sensitivity for safety assessment}

We consider a system with resistance $R$ and load $S$ and two model uncertainties $X_R$ and $X_S$. The system fails when $X_S S$ exceeds $X_R R$:
\begin{equation}\label{eq:Failure_example1}
	F=\left\{X_R R \le X_S S\right\}.
\end{equation}


A limit-state function describing this failure event is
\begin{equation}
g(\mathbf{X})=\ln X_R + \ln R - \ln X_S -\ln S.
\end{equation}

The probabilistic model of $\mathbf{X}=[R;S;X_R;X_S]$ is summarized in Table \ref{tab:prob_model_example1}. 

Because all random variables are modelled as lognormal, the probability of failure can be computed analytically.
It is 
\begin{equation} \label{eq:p_F_problem1}
p_F=\Phi \left[-\frac{\mu_{\ln X_R} + \mu_{\ln R} - \mu_{\ln X_S} -\mu_{\ln S}}{\sqrt{\sum_{i,j}\mathbf{C}_{ij}}}\right],
\end{equation}
wherein $\mathbf{C}$ is the covariance matrix of $\ln\mathbf{X}$.

The conditional probabilities of failure $p_F(x_i)$ are evaluated by first computing (analytically) the conditional moments of $\ln \mathbf{X}_{-i}$ given $X_i=x_i$. $p_F(x_i)$ is then obtained by analogy with Eq. \ref{eq:p_F_problem1}.

\begin{table}[]
	\centering
	\caption{Probabilistic model of example 1.}
	\label{tab:prob_model_example1}
	\small{
	\begin{tabular}{lllll}
		\hline
		 Parameter & Distribution & Mean & c.o.v. \\ \hline
		Resistance $R$          & lognormal       &    $100$  &  $0.2$                    \\
		Load $S$         & lognormal       &    $40$  & $0.25$               \\ 
		Model uncertainty $X_R$  & lognormal& $1$ & $0.1$
		\\ 
		Model uncertainty $X_S$  & lognormal& $1$ & $0.2$
		\\ \hline
		&              &      &                    &                                
	\end{tabular}
}
\end{table}

The EVPPI is evaluated (a) assuming independent input random variables and (b) assuming that the input random variables follow a Gaussian copula, with correlation matrix of $\mathbf{X}$ equal to
$$
\mathbf{R}_{\mathbf{X} \mathbf{X}} = 
\begin{bmatrix}
1 & 0 & 0.5 & 0 \\
0 & 1 & 0 & 0.5 \\
0.5 & 0 & 1 & 0.5 \\
0 & 0.5 & 0.5 & 1 
\end{bmatrix}.
$$
$\mathbf{C}$, the correlation matrix of $\ln \mathbf{X}$, can be evaluated in closed form in terms of $\mathbf{R}_{\mathbf{X} \mathbf{X}}$ \cite{Liu1986,Ditlevsen1996a}.

The probability of failure in the case of independent input random variables is $\Pr(F)=7.4\cdot10^{-3}$ and in the case of dependent input random variables it is $\Pr(F)=1.7\cdot10^{-2}$.

Cost values are $c_F=10^8$ and either $c_r=10^5$ or $c_r=10^6$, corresponding to cost ratios of either $\frac{c_r}{c_F}=10^{-3}$ or $\frac{c_r}{c_F}=10^{-2}$.

Figure \ref{fig:ConditionalPDF_Example_1} shows the relation between the prior PDF $f_R$, the conditional PDF $f_{R|F}$ and the resulting conditional failure probability $p_F(r)$ for the input random variable $R$ (independent case). Figure \ref{fig:CVPPI_linear_Example} shows the resulting $CVPPI(r)$ for $\frac{c_r}{c_F}=10^{-2}$.


%

\subsubsection{VoI sensitivity for reliability-based design}

For design optimization, we consider a modified version of the limit-state function,
\begin{equation}
g(\mathbf{X})=a X_R R - X_S S,
\end{equation}
in which $a$ is the design parameter affecting the overall resistance of the system. 

The probability of failure in function of $a$ is
\begin{equation}
p_F(a)=\Phi \left[-\frac{\ln a + \mu_{\ln X_R} + \mu_{\ln R} - \mu_{\ln X_S} -\mu_{\ln S}}{\sqrt{\sum_{i,j}\mathbf{C}_{ij}}}\right].
\end{equation}

For simplicity, we assume that the design cost is a linear function of $a$:
\begin{equation}
c_d(a) = c_\delta a.
\end{equation}
Two cases are evaluated, with cost factor equal to $c_\delta=10^5$ and $c_\delta=10^6$.

The optimal design a-priori is found by Eq. \ref{eq:prior_optimization_design} and reported in Table \ref{tab:optimal_design_example_1}, together with the associated failure probability.

\begin{table}[]
	\centering
	\caption{Optimal designs, together with the associated failure probabilities.}
	\label{tab:Optimal_designs_exmple_1}
	\small{
		\begin{tabular}{lllll}
			\hline
			& \multicolumn{2}{c}{$c_d(a)=10^5a$} & \multicolumn{2}{c}{$c_d(a)=10^6a$} \\
			& independent $\mathbf{X}$ & dependent $\mathbf{X}$ & independent $\mathbf{X}$ & dependent $\mathbf{X}$ \\ \hline
			$a_{opt}$ & $1.57$ & $1.87$  & $1.23$ &  $1.41$       \\
			$p_F(a_{opt})$         & $2.7\cdot10^{-4}$  & $2.2\cdot10^{-4}$  & $1.5\cdot10^{-3}$  & $2.0\cdot10^{-3}$                
			\\ \hline                        
		\end{tabular}
	}
	\label{tab:optimal_design_example_1}
\end{table}

%
%

\subsubsection{Results and discussion}

%

Tables \ref{tab:VPPI_example1_cost_model2} and \ref{tab:VPPI_example1_dependent} summarize the resulting decision sensitivities for selected cost models. In addition, FORM sensitivities are reported in terms of the squared $\alpha$-factors or the $\gamma^2$, which are a generalization of the $\alpha$-factors for dependent inputs \cite{kiureghian2005,Kim2018b}. The EVPPI are stated as normalized values
to facilitate a comparison with the squared $\alpha$-factors.

\begin{table}[]
	\centering
	\caption{Normalized sensitivity metrics for the component example with independent input random variables.}
	\label{tab:VPPI_example1_cost_model2}
	\small{
		\begin{tabular}{lcccccc}
			\hline
			&\multicolumn{4}{c}{Normalized EVPPI} & FORM    \\ \cline{2-5} 
			&safety & safety  &  design  & design & $\alpha_i^2$ \\
			&$\frac{c_r}{c_F}=10^{-3}$ & $\frac{c_r}{c_F}=10^{-2}$ &  $\frac{c_\delta}{c_F}=10^{-3}$  & $\frac{c_\delta}{c_F}=10^{-2}$  \\ \hline
			Resistance $R$  &  $ 25\% $  & $ 27\% $ & $ 26\% $  & $ 26\% $ & $ 26\% $  &   \\
			Load $S$         & $ 49\% $  & $ 35\% $ & $ 41\% $   & $ 41\% $ & $ 41\% $       \\ 
			Model uncertainty $X_R$  &$ 0.5\% $   &  $ 10\% $ & $ 7\% $	& $ 6\% $   & $ 7\% $  	\\ 
			Model uncertainty $X_S$  & $ 25\% $   &  $27\% $ & $ 26\% $ & $ 26\% $   & $ 26\% $ 
			\\ \hline                        
		\end{tabular}
	}
\end{table}

\begin{table}[]
	\centering
	\caption{Normalized sensitivity metrics for the component example with dependent input random variables.}
	\label{tab:VPPI_example1_dependent}
	\small{
		\begin{tabular}{lcccccc}
			\hline
			&\multicolumn{4}{c}{Normalized EVPPI} & FORM    \\ \cline{2-5} 
			&safety & safety  &  design  & design & $\gamma_i^2$   \\
			&$\frac{c_r}{c_F}=10^{-3}$ & $\frac{c_r}{c_F}=10^{-2}$ &  $\frac{c_\delta}{c_F}=10^{-3}$  & $\frac{c_\delta}{c_F}=10^{-2}$  \\ \hline
			Resistance $R$  &  $ 15\% $  & $ 26\% $ & $ 23\% $  & $ 23\% $ & $ 26\% $  &   \\
			Load $S$         & $ 61\% $  & $ 41\% $ & $ 46\% $   & $ 46\% $ & $ 41\% $       \\ 
			Model uncertainty $X_R$  &$ 0.0\% $   &  $ 4\% $ & $ 4\% $	& $ 4\% $   & $ 7\% $  	\\ 
			Model uncertainty $X_S$  & $ 24\% $   &  $29\% $ & $ 28\% $ & $ 28\% $   & $ 26\% $ 
			\\ \hline                        
		\end{tabular}
			}
\end{table}

As expected for this linear Gaussian case, the ranking coincides for all EVPPI cases and the FORM sensitivities. Note that in the independent case, this equivalence follows from Sections \ref{sec:FORM_approx} and \ref{sec:FORM_design_case}.

However, for the safety decision, the relative contribution varies and is notably different for cost ratio $\frac{c_r}{c_F}=10^{-3}$. In this case, the EVPPI of input $X_R$ is negligible, which is due to the fact that at this cost ratio, it is very unlikely that learning about $X_R$ will change the a-priori optimal decision. 

The design EVPPI results of the independent case are similar to the squared $\alpha$-factors, which confirms the conclusions in Section \ref{sec:FORM_design_case}. 
Generally, the results for the EVPPI sensitivity of the reliability-based design case are independent of the ratio $\frac{c_\delta}{c_F}$ for this linear Gaussian model. This is in agreement with the findings reported in the Annex. 

\begin{figure}
	\centering	\makebox[0pt]{\includegraphics[width=170mm]{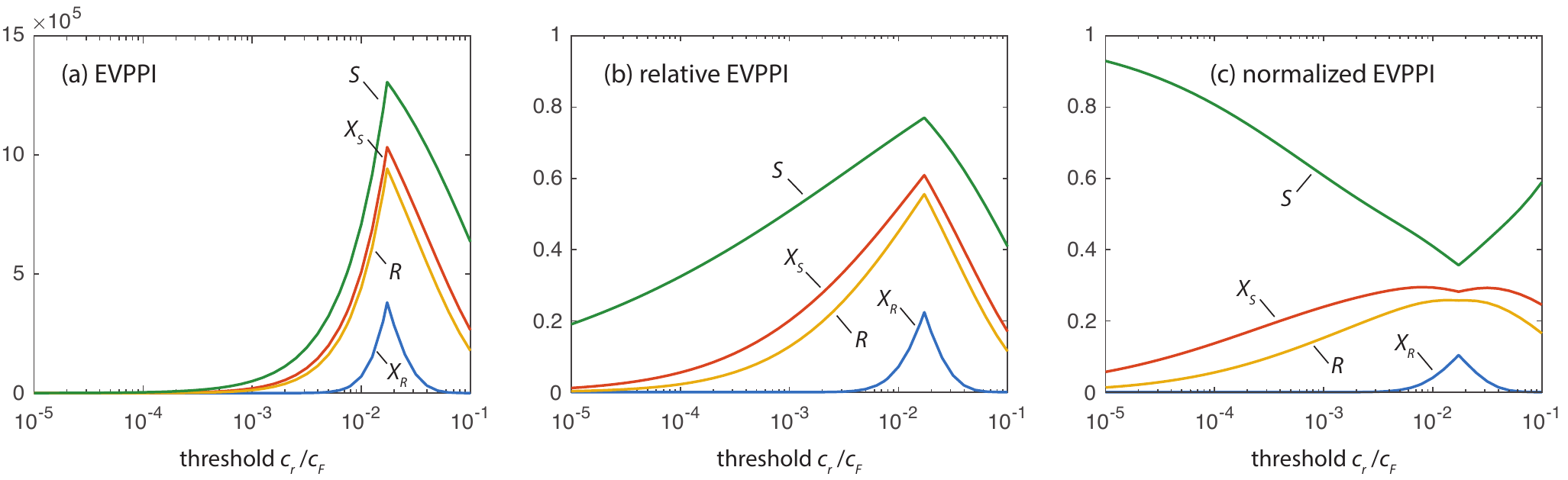}}
	\caption{EVPPI of the safety decision case in function of the threshold $\frac{c_r}{c_F}$ (case with dependent inputs $\mathbf{X}$). The left panel shows the absolute value, the middle panel the relative EVPPI (normalized with the EVPI following Eqs. \ref{eq:relative_EVPPI} and \ref{eq:safety_EVPI}) and the right panel shows the EVPPI normalized to sum to one (Eq. \ref{eq:normalized_EVPPI}).}
	\label{fig:EVPPI_in_function_of_threshold}
\end{figure}

Figure \ref{fig:EVPPI_in_function_of_threshold} plots the EVPPI of the safety decision case in function of $\frac{c_r}{c_F}$ for the case with dependent inputs. (The plot for independent inputs looks similar.) It compares the absolute value of the EVPPI with the relative and the normalized EVPPI. Both the absolute and the relative EVPPI are highest with a decision threshold $\frac{c_r}{c_F}$ close to the probability of failure $p_F$. This is not surprising, since in this case it is more likely that learning an input parameter will change the optimal decision. The normalized EVPPI shows that when $\frac{c_r}{c_F}$ is close to $p_F$, the differences between the EVPPIs are smaller. Near the threshold, even learning the less important input random variables still has considerable value. The further $\frac{c_r}{c_F}$ deviates from $p_F$, the larger the differences between EVPPIs. In the extreme, the normalized EVPPI of the most important input random variable will approach $1$. 




\subsubsection{Monte Carlo approximation for the safety decision case}
To investigate the accuracy of the Monte Carlo approach of Section \ref{sec:failure_samples}, we determine the EVPPI with $10^2$ and $10^3$ failure samples. This is the typical range of the number of samples in the failure domain one would expect from a sampling-based reliability evaluation.   

The results are summarized in Table \ref{tab:VPPI_example1_sample_evaluation}. Each analysis is repeated $100$ times, to obtain an approximate mean and standard deviation of the MCS estimate. The results show that for this example the accuracy of the estimate is sufficient for as low as $n_F=100$ independent failure samples. 

\begin{table}[]
	\centering
	\caption{EVPPI for the safety-assessment case with $\frac{c_r}{c_F}=10^{-2}$, $c_F=10^8$ and independent input random variables. Exact value and sample-based value with $n_F$ failure samples. For the latter, the mean value and coefficient of variation (in parenthesis) are evaluated by $100$ repetitions of the analysis.}
	\label{tab:VPPI_example1_sample_evaluation}
	\begin{tabular}{lcccc}
		\hline
		&\multicolumn{3}{c} {EVPPI} &  \\ \cline{2-4}
		&exact & $n_F=10^3$  &  $n_F=10^2$  \\ \hline
		Resistance $R$  & $ 349 $  & $ 344 (3.2\%) $ & $ 337(7.1\%) $   \\
		Load $S$         & $ 454 $   & $ 443 (2.3\%) $ & $ 443 (6.0\%)$   \\ 
		Model uncertainty $X_R$ & $ 131  $	& $ 131 (8.3\%) $   & $ 127 (22.4\%) $  	\\ 
		Model uncertainty $X_S$  &  $ 349 $ & $ 344 (3.2\%)$   & $ 339 (7.9\%) $  
		\\ \hline                        
	\end{tabular}
\end{table}

\subsubsection{Discrete approximation for the design decision case} \label{sec:example_discrete_approx}
We investigate the discrete approximation of design choices following Section \ref{sec:discretized_design}. Results are shown here for the setting $c_\delta=10^5$ and independent input random variables.

Discrete design parameter choices between $a_1=0.5$ and $a_m=2.0$ are considered, with step size $\Delta a=\frac{1.5}{m-1}$. The approximation of the posterior expected loss in case of $m=4$ is shown in Figure \ref{fig:Posterior_loss_discrete_Ex1}.

\begin{figure}
\centering	\makebox[0pt]{\includegraphics[width=80mm]{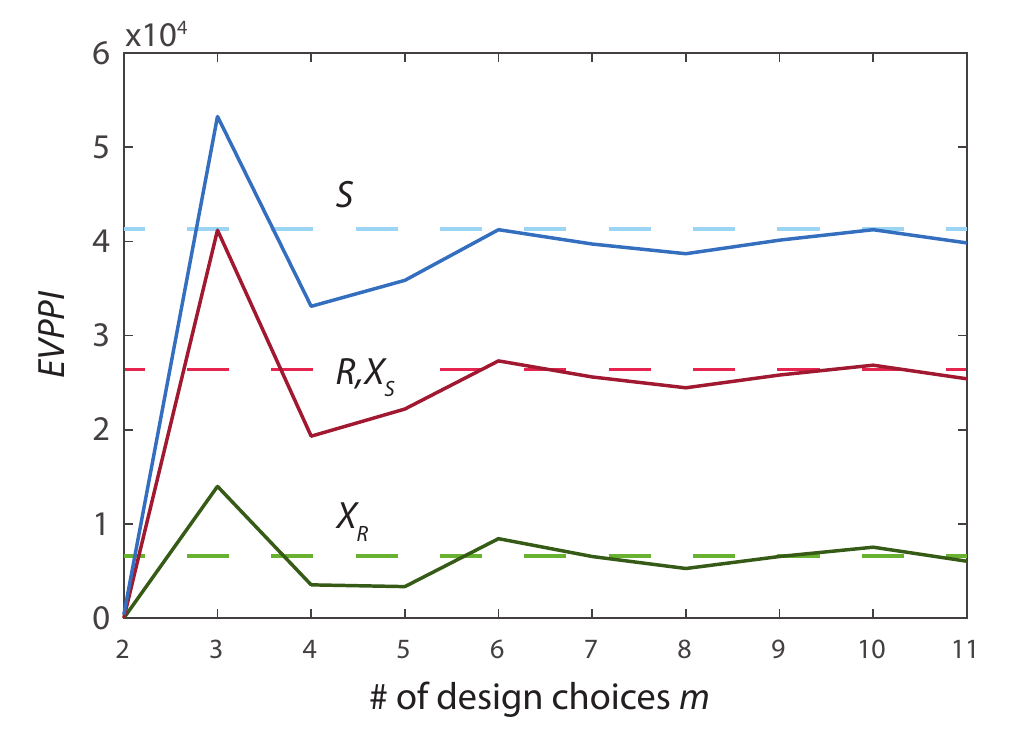}}
\caption{Convergence of the EVPPI with number of discrete values $m$ of the design choice. (EVPPI evaluated for the reliability-based design case with $c_\delta=10^5$ and independent input random variables.) }
\label{fig:Discrete_approx_error}
\end{figure}

Figure \ref{fig:Discrete_approx_error} shows the convergence of estimated EVPPI to the exact values with increasing $m$. While an accurate evaluation of the absolute EVPPI here requires $m\ge 6$ discrete design choices, the normalized EVPPI is already evaluated reasonably well with $m$ as low as $3$, see Table \ref{tab:VPPI_example1_discretization_evaluation}.

\begin{table}[]
	\centering
	\caption{Normalized EVPPI for the reliability-based design case with $c_\delta=10^{5}$ and independent input random variables. Exact value and approximation with a discrete representation of the design choice with $m$ values.}
	\label{tab:VPPI_example1_discretization_evaluation}
		\begin{tabular}{lcccccc}
			\hline
			&\multicolumn{4}{c} {Normalized EVPPI} &  \\ \cline{2-5}
			&exact & m=3  &  m=4  & m=6  \\ \hline
			Resistance $R$  & $ 26\% $  & $ 27\% $ & $ 26\% $  &  $26\%$ \\
			Load $S$         & $ 41\% $   & $ 36\% $ & $ 44\% $  &$39\%$     \\ 
			Model uncertaint $X_R$ & $ 7\% $	& $ 9\% $   & $ 5\% $ &$8\%$ 	\\ 
			Model uncertainty $X_S$  &  $ 26\% $ & $ 27\% $   & $ 26\% $ &$26\%$ 
			\\ \hline                        
		\end{tabular}
\end{table}

\subsection{Example 2: Non-linear limit state} \label{sec:nonlinear_example}
We consider a short column subjected to biaxial bending moments $M_1$ and $M_2$ and axial force $P$ that was previously investigated in \cite{kiureghian2005,Kim2018b}. The failure is defined by the limit-state function:
\begin{equation}
	g(\mathbf{X})=1-\frac{M_1}{s_1 Y}-\frac{M_2}{s_2 Y}-\left(\frac{P}{A Y}\right)^2.
\end{equation}
The random variables are $\mathbf{X}=[M_1;M_2;P;Y]$, with $Y$ being the yield strength of the material. Their joint probability distribution is a Gaussian copula model as in \cite{Kim2018b} with the parameters given in Table \ref{tab:prob_model_example2}. The deterministic parameters are the flexural moduli of the plastic column section $s_1=0.03\text{m}^3$ and $s_2=0.015\text{m}^3$ and the column cross section $A=0.190\text{m}^2$. 

\begin{table}[]
	\centering
	\caption{Probabilistic model of example 2, including correlation matrix $\mathbf{R_{XX}}$.}
	\label{tab:prob_model_example2}
	\small{
		\begin{tabular}{llllcccc}
			\hline
			Parameter & Distribution & Mean & c.o.v. & \multicolumn{4}{c} {$\mathbf{R_{XX}}$} \\ \cline{5-8}  &&&&$M_1$&$M_2$&$P$&$Y$\\\hline
			Bending moment $M_1$ [kNm]      & normal       &    $250$  &  $0.3$  &$1.0$ &$0.5$ & $0.3$ & $0.0$                  \\
			Bending moment $M_2$ [kNm]    & normal       &    $125$  & $0.3$  &$0.5$ &$1.0$ & $0.3$ & $0.0$                        \\ 
			Axial force $P$ [kN]  & Gumbel& $2500$ & $0.2$ &$0.3$ &$0.3$ & $1.3$ & $0.0$           
			\\ 
			Yield strength $Y$ [N/$\text{mm}^2$]  & Weibull& $40$ & $0.1$ &$0.0$ &$0.0$ & $0.0$ & $1.0$           
			\\ \hline
			&              &      &                    &                                
		\end{tabular}
	}
\end{table}

The reliability is evaluated with crude Monte Carlo with $10^6$ samples. The MC estimate of the probability of failure is $0.0094$ with $95\%$ credible interval $[0.0092,0.0096]$. 

We consider only the safety assessment case. The EVPPI is determined based on the MC samples in the failure domain, following Section \ref{sec:failure_samples}. Figure \ref{fig:KDE_fit_illustration} exemplary shows the KDE fit of the conditional distribution of $M_1$ given failure. 

Different cost ratios $\frac{c_r}{c_F}$ are considered to reflect different decision situations. The results are summarized in Table \ref{tab:results_example2}, together with results from the sensitivity metrics evaluated in \cite{Kim2018b}.

\begin{table}[]
	\centering
	\caption{Normalized sensitivity indices for example 2. FORM and GRIM results are taken from \cite{Kim2018b} and are provided for comparison.}
	\label{tab:results_example2}
	\small{
		\begin{tabular}{llllll}
			\hline Parameter
			& Nor. EVPPI & Nor. EVPPI & Nor. EVPPI & FORM & GRIM  \\ & ($\frac{c_r}{c_F} = 10^{-3}$) & ($\frac{c_r}{c_F} = 10^{-2}$) & ($\frac{c_r}{c_F} = 10^{-1}$) & $\alpha^2$  &  $\gamma^2$  \\ \hline
			$M_1$  & $21\%$       &    $22\%$  &$9\%$ &  $7\%$ &$8\%$                \\
			$M_2$  & $23\%$       &    $23\%$  & $9\%$ & $7\%$ & $6\%$               \\ 
			$P$ & $8\%$& $24\%$ &$28\%$& $20\%$ & $25\%$
			\\ 
			$Y$  & $48\%$& $31\%$ & $54\%$ & $65\%$ & $61\%$
			\\ \hline
			&              &      &                    &                                
		\end{tabular}
	}
\end{table}

The resulting EVPPI indicate that the decision sensitivities can vary strongly among different decision situations. For the situation $\frac{c_r}{c_F} = 10^{-1}$, the EVPPI is more or less comparable to the sensitivies obtained with the indices from FORM and GRIM (Generalized Reliability Importance Measures). This corresponds to the case where the cost of repair is large and hence the a-priori decision is not to repair. $\frac{c_r}{c_F} = 10^{-2}$ is very close to the a-priori probability of failure of $0.0094$. In this case, learning the value of any of the four inputs is likely to change the decision from not repairing to repairing, which explains the more uniform sensitivity values. In the case $\frac{c_r}{c_F} = 10^{-3}$, the EVPPI leads to a different ranking. In this case, the prior decision is to repair. The probability of changing this decision depends on the other tails of the distributions than in cases $\frac{c_r}{c_F} = 10^{-1}$ and $\frac{c_r}{c_F} = 10^{-2}$. This explains the differences between these two decision situations, and also highlights that the EVPPI does account for the full distribution


%
%
%
%
%

\subsection{Example 3: Foundation reliability} \label{sec:foundation_example}
\begin{figure}
	\centering
	\def\svgwidth{.5\textwidth}
	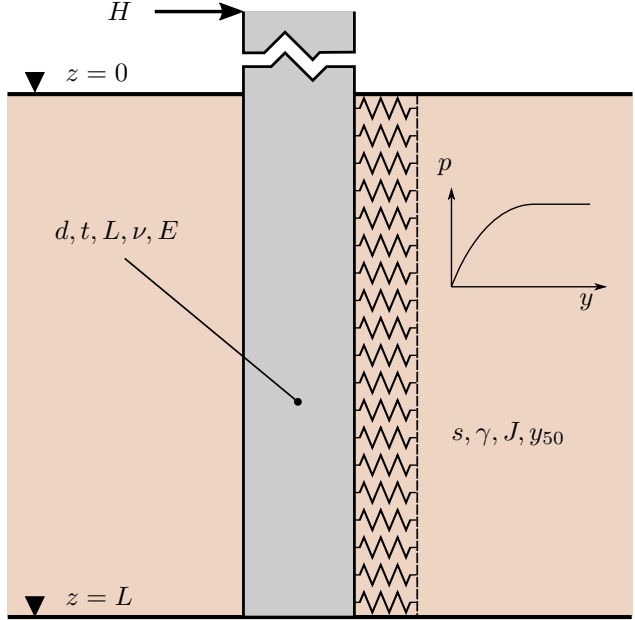{}
	\caption{Wind turbine monopile foundation  \cite{Depina2017}.}
	\label{fig:monopile}
\end{figure}

We consider a monopile foundation of an offshore wind turbine in stiff, plastic soil (Figure \ref{fig:monopile}). The soil-structure interaction is described by a nonlinear finite element model.
Details regarding the engineering model setup can be found in \cite{Depina2017}.
Deterministic parameters of the monopile foundation are its depth $L = 30$ m, diameter $d = 6$ m, wall thickness $t=0.07$ m, Poisson ratio $\nu = 0.3$ and Young's modulus $E = 2.1\cdot10^{5}$ MPa.
The uncertain inputs comprise the Gumbel-distributed lateral load $H$ and its distribution parameters $a_H$ and $b_H$ as well as the undrained shear strength $s$ and its hyperparameters $s_0$ and $s_1$. $s$ is modelled as a 1-D, non-stationary, log-normal random field with a linear drift along the $z$-coordinate. It is discretized using the midpoint method \cite{der1988stochastic} and an 82-dimensional standard-normal random vector $\boldsymbol{\xi}$. Its hyperparameters are the ground-level shear strength $s_0$ (at $z=0$) and the mean shear strength slope $s_1$.
The latter two can be thought of as representing inter-site variability while $\boldsymbol{\xi}$ models intra-site variability.
The probabilistic inputs are summarized in Table \ref{tab:monopile}. A detailed description of the probabilistic model setup is provided in \cite{Ehre2020,Jiang2018}.

\begin{table}[!ht]
	\caption{ Input variable definitions of the monopile foundation.}
	\label{tab:monopile}
	\centering
	\begin{tabular}{p{2cm}ccc}
		\hline
		Input & Distribution & Mean $\mu$  & CoV $\delta$\\
		\hline
		$\boldsymbol{\xi}$ [-]& Standard-Normal & $\boldsymbol{0}$ & $n.d.$ ($\boldsymbol{\Sigma}_{\boldsymbol{\xi}\boldsymbol{\xi}} = \boldsymbol{I}_{82 \times 82}$)\\
		$s_0$ [kPa]& Log-Normal & $33.7094$ & $0.3692$\\
		$s_1$ [kPa]& Log-Normal & $0.7274$ & $0.8019$\\
		$H$ [kN]& Gumbel & $\mu_{H|a_H,b_H}$ & $\delta_{H|a_H,b_H}$\\
		$a_H$ [kN] & Log-Normal & $2274.97$ & $0.2$ \\
		$b_H$ [kN]& Log-Normal &  $225.02$&  $0.2$\\
		\hline
	\end{tabular}
\end{table}

Failure occurs when the maximum stress in the foundation, $\sigma_{\text{max}}(\mathbf{X})$ exceeds $\sigma_{\text{crit}} =100~\text{MPa}$. 
The corresponding limit-state function is
	\begin{equation} 
	\label{e:lsf_monopile}
	g(\mathbf{X}) = \sigma_{\text{crit}} - \sigma_{\text{max}}(\mathbf{X}).
	\end{equation}

The analysis is performed with subset simulation with $10^4$ samples per level \cite{Au2001,papaioannou2015mcmc}. To evaluate the coefficient of variation of the subset simulation estimate, we repeat the computation $100$ times.
The a-priori failure probability is estimated as $P(F) \approx 3.66 \cdot 10^{-4}$ with a c.o.v. of $0.08$. This result is also verified with crude Monte Carlo based on $2.8\cdot10^6$ samples, which results in a $95\%$ credible interval $[3.33,3.78] \cdot 10^{-4}$.


\subsubsection{VoI sensitivity for safety assessment}
The safety assessment scenario is evaluated for a cost of failure $c_F = 10^8$ and two values of the cost of repair, $c_r = [10^4;10^5]$. Hence the resulting cost ratios are $\frac{c_r}{c_F} = [10^{-3};10^{-4}]$.


The EVPPI is estimated with failure samples according to Section \ref{sec:failure_samples} for all input random variables except $\boldsymbol{\xi}$ and $H$. The decision sensitivity of the former cannot be estimated reliably with the proposed method, due to its $82$ dimensions, but we know from \cite{Ehre2020} that its importance is small. The decision sensitivity of the latter can be estimated, but is not meaningful, since it is not possible to reduce uncertainty on $H$ directly. It is only possible to reduce the uncertainty of its distribution parameters, $a_H$ and $b_H$. 

The conditional density in the numerator of Eq. \ref{eq:pF_cond_Bayes} is estimated with the $10^4$ samples that are produced in the final level of subset simulation. Note that these samples are dependent, because of the MCMC algorithm that is used to propagate samples from one level to the next in subset simulation. 
The results are reported in Table \ref{tab:VPPI_example5}. In addition to the normalized EVPPI, also the reliability-oriented Sobol' indices of ${\log p_F}$ introduced in \cite{Ehre2020} are presented for comparison. 







\begin{table}[]
	\centering
	\caption{Normalized sensitivity measures for the monopile. In brackets: coefficient of variation of the estimate (first position), absolute estimated value (second position). The coefficients of variation of the safety cases are estimated based on 100 repeated runs of subset simulation.}
	\label{tab:VPPI_example5}
	\small{
		\begin{tabular}{lccccc}
			\hline
			&\multicolumn{4}{c} {Normalized EVPPI} &  Sobol' \cite{Ehre2020}\\
			&safety   & safety  &  design  & design &   \\
			& $\frac{c_r}{c_F}=10^{-4}$   & $\frac{c_r}{c_F}=10^{-3}$   &  linear $c_d(a)$ & quadratic $c_d(a)$ &  \\ \hline
			$s_0$               & $0\%~(1.6 ; 4.7)$
			&   $2\%~(0.84; 542)$ & $ 8\%~(-; 4971)$   & $ 8\%~(-; 2790) $ & $ 0.58\% $  \\ 
			$s_1$               & $0\%~(1.81; 2.5)$ &    $34\%~(0.23; 7414)$ & $ 19\%~(-; 11434) $
			& $ 14\%~(-; 5084) $  & $  4.59\% $ 
			\\ 
			$a_H$               & 
			$36\%~(0.63; 678)$   & $36\%~(0.16; 8019)$ & $ 41\%~(-; 24351) $ & $ 41\%~(-; 14542) $  & $ 22.70\% $ 
			\\
			$b_H$               & $63\%~(0.31; 1176)$   &  $28\%~(0.26; 6071)$ & $ 32\%~(-; 18792) $ & $ 36\%~(-; 12662) $ & $ 70.39\% $ 
			\\ \hline                        
		\end{tabular}
	}
\end{table}
In the case $\frac{c_r}{c_F}=10^{-4}$, the EVPPI-based ranking is identical to the one of obtained with the reliability-oriented Sobol' indices. In the case $\frac{c_r}{c_F}=10^{-3}$, the EVPPI-based ranking differs. In particular, $s_1$ becomes considerably more relevant as the repair cost is increased and $b_h$ becomes less important.

\subsubsection{VoI sensitivity for reliability-based design}
\begin{figure}
	\centering	\makebox[0pt]{\includegraphics[width=50mm]{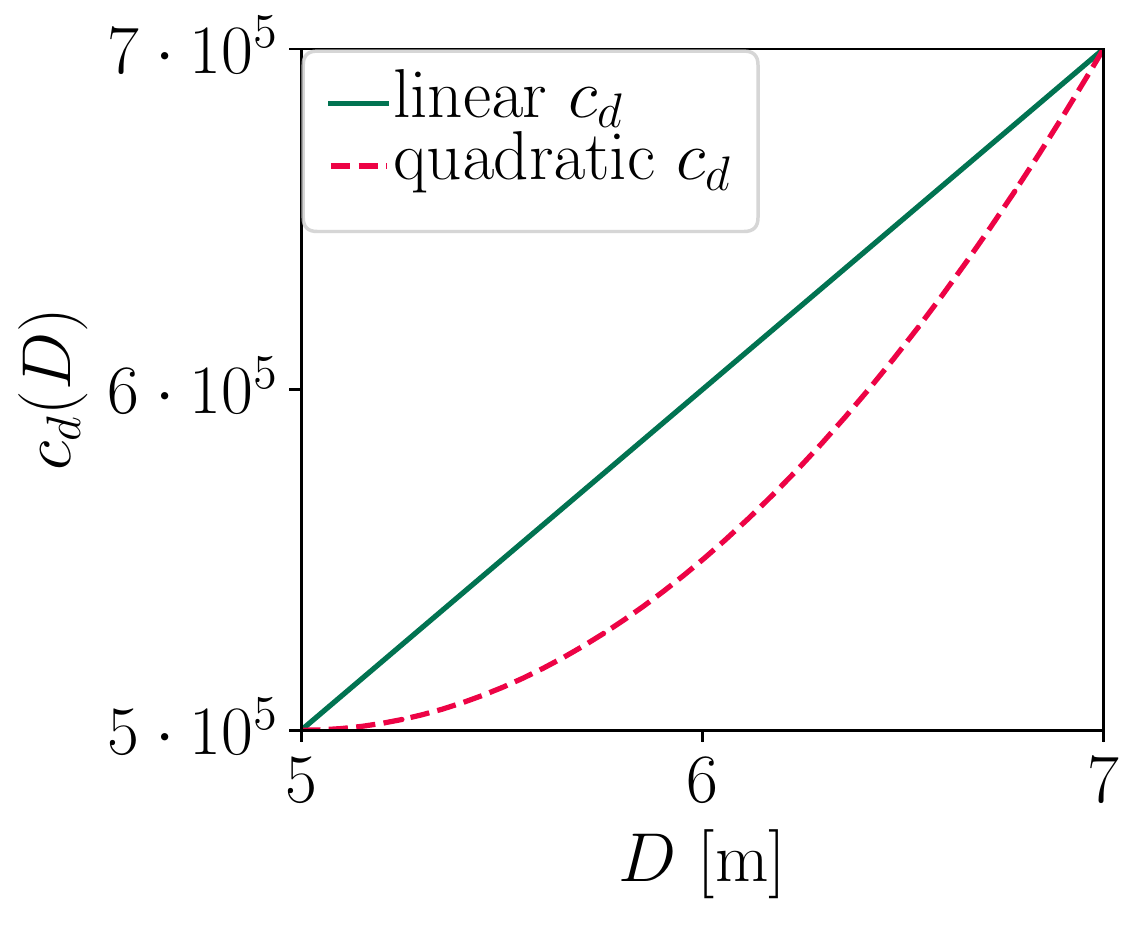}}
	\caption{Linear and quadratic design cost models. 
	}
	\label{fig:cost_models_ex5}
\end{figure}

We consider the monopile diameter $d \in \mathbb{D}$ as the design parameter, with the feasible domain $\mathbb{D} = [5,7]$ m. The cost of failure is $c_F = 10^8$ and the design cost is modeled with either a linear cost model $c_d(d) = 10^5\cdot d$ or a quadratic cost model
$c_d(d) = (5 + 2 \cdot ((d-5)/2)^2) \cdot 10^5$, such that for both models, $5 \cdot 10^5 \le c_d \le 7 \cdot 10^5$ on $\mathbb{D}$. The two models are plotted in Figure \ref{fig:cost_models_ex5}.
We discretize the feasible domain for the design parameter in $m=201$ steps of $\Delta d = 0.01m$ as $\{5.00,5.01,5.02,\dots,7.00\}$m and run $201$ subset simulations to obtain $10^4$ failure samples associated with each discrete design. Based on these failure samples, we compute the conditional failure probabilities associated with all variables but $\boldsymbol{\xi}$ and $H$ according to Eq.~\ref{eq:pF_cond_Bayes} for each discrete design.
The EVPPI results for the reliability-based design case are summarized in Table \ref{tab:VPPI_example5}.
The linear cost model produces more expensive designs on average thus yielding larger absolute EVPPI than the quadratic model. 
The difference in normalized EVPPI (percentages) between the two cost models is negligible. 



\section{Discussion}
Many sensitivity metrics exist in the literature, also in the context of reliability assessment. The interpretation of these metrics is often difficult, even for experts \cite{Saltelli2019}. In this sense, the decision-theoretic sensitivity metrics have a clear advantage: They have a straightforward interpretation in the context of a specific decision situation. This makes them applicable to all models and applications, including system reliability problems and problems with dependent inputs. 

In this paper, we have shown how the expected value of partial perfect information (EVPPI) can be evaluated by post-processing the results of a reliability analysis. For the safety decision case, no additional model evaluations are necessary, hence the computation of EVPPI is cheap. For the design decision case, additional model evaluations are necessary, so the metric might not be convenient when models are computationally expensive. Even if not sufficient computational ressources are available to evaluate the EVPPI for the design decision, the consideration of the decision-theoretic basis of the EVPPI metric can nevertheless be helpful for interpreting other sensitivity metrics. 

If suitable surrogate models exist, efficient computation of the EVPPI in the design decision case might be possible. We have experimented with such a surrogate model approach, based on \cite{Ehre2020}, and have applied it to the example of Section \ref{sec:foundation_example}. However, we have found that the computation of the EVPPI requires a surrogate with high accuracy over a wide range of the input random variable's outcome space. Standard surrogates used for reliability analysis do not have this property. The development of suitable surrogates is left for future research.

Our paper focuses on computational strategies for evaluating the EVPPI and provides only limited guidance on how to select the decision case and the associated cost/threshold parameters in specific applications. In some projects the analyst faces a situation closely resembling the decision cases presented here, in other projects less so. In the following, we discuss two applications with different decision situations, to indicate how the presented decision cases can nevertheless be utilized as a proxy if the cost/threshold parameters are chosen appropriately. 

Firstly, in the analysis presented here we have assumed that the decision is taken based on minimizing expected loss. However, regulatory constraints often result in a decision that is based on a reliability criterion in the form of a safety threshold $p_F^T$. In this case, one aims at minimizing the expected life-cycle costs under the constraint $p_F\le p_F^T$. This case is also considered in \cite{Borgonovo2017}, who suggest to find a utility function such that the optimal decision is repair ($a_r$) if $\Pr(F)>p_F^T$ and do nothing ($a_0$) otherwise. In our framework, this translates to setting $\frac{c_r}{c_F}=p_F^T$ and then applying the methods of Section \ref{sec:Safety_assessment}.

If the costs $c_r$ and $c_F$ are known, setting $\frac{c_r}{c_F}=p_F^T$ might not reflect well the decision situation. Instead, one could include the constraint $p_F^T$ into the decision analysis. This might, however, have the effect that the EVPPI can become negative, e.g., if the a-priori probability of failure is slightly below the acceptable $p_F^T$. In this case, the decision maker might prefer not to gain additional information, as that might change the probability of failure estimate to an unacceptable level and incur the repair costs $c_r$. Such information avoidance has been discussed in the literature \cite{Pozzi2020}. For the purpose of sensitivity analysis, we suggest to avoid the issue, by either ignoring the safety criterion $p_F^T$ in the EVPPI calculation, or by setting $\frac{c_r}{c_F}=p_F^T$. 

Secondly, scientific analyses are often undertaken without a specific decision in mind \cite{borgonovo2021probabilistic}. Or the analyst might find it difficult to specify a loss function for other reasons. In these cases, an EVPPI based on the safety decision case might be computed, with varying threshold $\frac{c_r}{c_F}$, as depicted in Figure \ref{fig:EVPPI_in_function_of_threshold}. Reporting sensitivities for varying $\frac{c_r}{c_F}$ has the advantage that the dependence of the sensitivity on the decision context is clearly communicated.

\section{Conclusion}
We demonstrate that the expected value of partial perfect information (EVPPI) is an informative sensitivity metric for factor prioritization in reliability assessments, which -- in contrast to other sensitivity metrics -- is straightforward to interpret. It is a measure of decision sensitivity and measures how learning an input random variable improves the decision taken based on the reliability assessment.
With the computational strategies proposed in this paper, the EVPPI of the safety decision case can be evaluated following a reliability analysis without any additional model runs. 

\section*{Acknowledgment}
We acknowledge support by the German Research Foundation (DFG) through Grants STR 1140/11-1 and PA 2901/1-1.


\bibliography{VOI_sensitivity}

\newpage
\section*{Annex: A FORM approximation for the reliability-based design case}

We postulate the following design limit-state function with design parameter $a$:
\begin{equation}\label{eq:design_LSF_FORM}
	g_d(\mathbf{X},a)=g(\mathbf{X})+a.
\end{equation}
We let $\beta_0$ denote the FORM reliability index for $a=0$ and $\boldsymbol{\alpha}=\left[\alpha_1,\dots,\alpha_n\right]$ the corresponding row vector of FORM sensitivities. The corresponding FORM approximation $G_1(\mathbf{U})$ of the limit-state function in standard normal space is
\begin{equation}
	G_1(\mathbf{U})=|| \nabla G(\mathbf{u}^*)|| \left(\beta_0-\boldsymbol{\alpha}\mathbf{U}\right). 
\end{equation}
$|| \nabla G(\mathbf{u}^*)||$ is the norm of the gradient of limit-state function $G$ in standard normal space, evaluated at the design point.
Without loss of generality, we set $|| \nabla G(\mathbf{u}^*)||=1$. 

We assume that $G_1(\mathbf{U})$ is also a reasonable approximation of $g(\mathbf{X})$ when $a\ne 0$, which is valid for mildly non-linear problems. Hence the approximation of $g_d(\mathbf{X},a)$ in standard normal space becomes 
\begin{equation}\label{eq:FORM_approx_design_LSF}
	G_d(\mathbf{U},a)=\beta_0-\boldsymbol{\alpha}\mathbf{U} + a.
\end{equation}

The corresponding approximation of the probability of failure in function of $a$ is
\begin{equation}
	p_F(a)\approx \Phi\left(-\beta_0-a\right).
\end{equation}

Additionally conditioning on $u_i$ results in (by analogy with Eq. \ref{eq:pF_ui})
\begin{equation}
	p_F(u_i,a)\approx \Phi\left(\frac{\alpha_i u_i-\beta_0-a}{\sqrt{1-\alpha_i^2}}\right).
\end{equation}

Using a linear cost function $c_d(a)=c_\delta a$, the optimal design choice following Eq. \ref{eq:prior_optimization_design} is
\begin{equation} 
	\begin{split}
		a_{opt}	&= \arg \min_a \, c_\delta a + p_F(a) c_F  \\
		& \approx \arg \min_a \, c_\delta a + \Phi\left(-\beta_0-a\right)  c_F.  
	\end{split}
\end{equation}

Without loss of generality, we choose $c_\delta$ such that $a_{opt}=0$ is the optimal design in the a-priori case, which results in
\begin{equation} \label{eq:FORM_design_cost}
	c_\delta=\varphi \left(-\beta_0\right)c_F,
\end{equation}
wherein $\varphi$ is the standard normal PDF.

The conditionally optimal design is 
\begin{equation} 
	\begin{split}
		a_{opt|U_i}(u_i) &= \arg \min_a \, c_\delta a + p_F(u_i,a) c_F  \\
		& \approx \arg \min_a \, c_\delta a + \Phi\left(\frac{\alpha_i u_i-\beta_0-a}{\sqrt{1-\alpha_i^2}}\right)   c_F \\ &= \alpha_i u_i-\beta_0+\sqrt{-2\ln\left(\sqrt{2\pi}\frac{c_\delta}{c_F}\sqrt{1-\alpha_i^2}\right)\left(1-\alpha_i^2\right)}.		  
	\end{split}
\end{equation}
Inserting Eq. \ref{eq:FORM_design_cost} results in
\begin{equation} 
		a_{opt|U_i}(u_i) =  \alpha_i u_i-\beta_0+\sqrt{\left[\beta_0^2-\ln\left(1-\alpha_i^2\right)\right]\left(1-\alpha_i^2\right)}. 
\end{equation}

For ease of notation, we introduce\begin{equation}
	b = \sqrt{\beta_0^2-\ln\left(1-\alpha_i^2\right)}.
\end{equation}

By inserting the above expressions into Eq. \ref{eq:VPPI_design}, the EVPPI follows as
\begin{equation} \label{eq:EVPPI_FORM_design}
		\begin{split}
		EVPPI_{X_i} &= c_\delta a_{opt} + p_F(a_{opt}) c_F   -  \text{E}_{U_i}\left[c_\delta a_{opt|U_i}(U_i) + p_F\{U_i,a_{opt|U_i}(U_i)\} c_F\right]  \\ 
		&=\Phi\left(-\beta_0\right)c_F- \text{E}_{U_i}  \left[ c_\delta  \left( \alpha_i u_i-\beta_0+b\sqrt{\left(1-\alpha_i^2\right)}	\right)
		+\Phi\left(-b\right) c_F \right] \\
		&=\left[\Phi\left(-\beta_0\right)+\varphi\left(-\beta_0\right)\left(\beta_0-b \sqrt{\left(1-\alpha_i^2\right)}\right) -\Phi\left(-b\right)
		\right]c_F,
	\end{split} 
\end{equation}
where we have utilized $a_{opt}=0$ and $\text{E}_{U_i}\left[U_i\right]=0$. 

\begin{figure}
	\centering	\makebox[0pt]{\includegraphics[width=140mm]{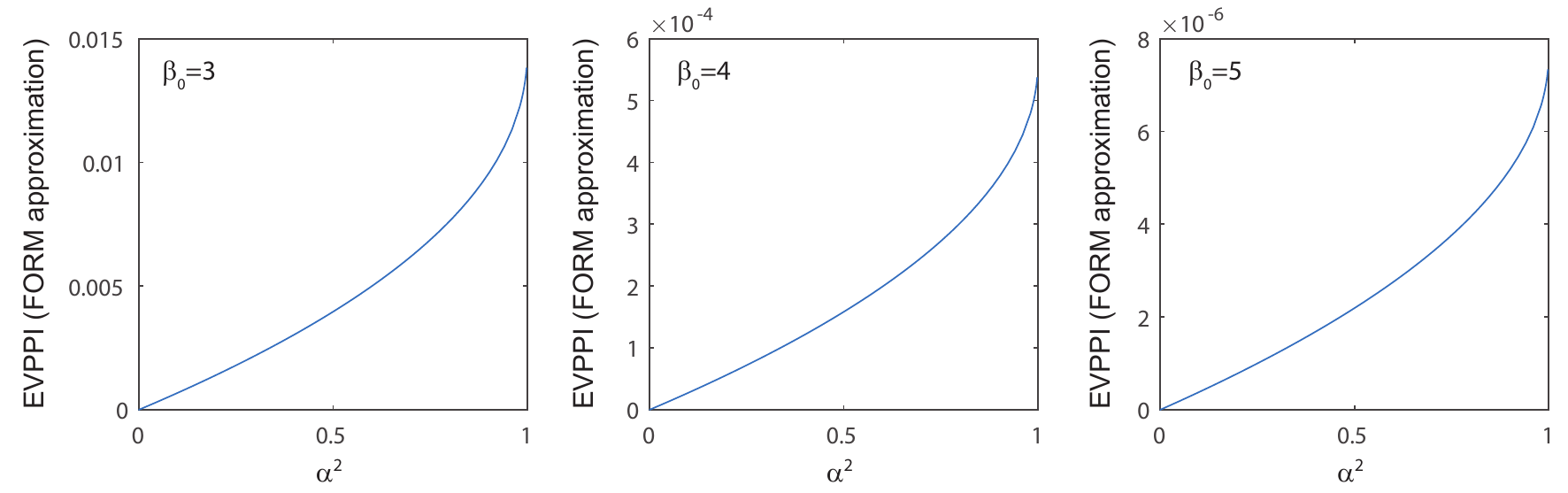}}
	\caption{EVPPI evaluated with the FORM approximation in function of the squared FORM $\alpha$-factors, for different levels of reliability $\beta_0$. The cost of failure is $c_F=1$ and $c_\delta$ is determined according to Eq. \ref{eq:FORM_design_cost}.}
	\label{fig:FORM_design_diffBeta}
\end{figure}

Figure \ref{fig:FORM_design_diffBeta} summarizes the resulting EVPPI in function of $\alpha_i^2$ for different $\beta_0$. The relative EVPPI, which is obtained by dividing the EVPPI with the value obtained at $\alpha^2=1$, is almost identical for different $\beta_0$, as shown in Figure \ref{fig:FORM_design_norm}. Note that the relative EVPPI does not depend on $c_F$ and $c_\delta$. Hence the relative EVPPI values shown in Figure \ref{fig:FORM_design_norm} are expected to be approximately valid for any reliability problem, which has a single design point and a limit-state function that is mildly non-linear in standard normal space.

\end{document}